\documentclass[11pt]{article}
\usepackage{latexsym,amssymb}
\usepackage{amsmath}
\usepackage{color}

\newtheorem{Theorem}{\bf Theorem}[section]
\newtheorem{Lemma}{\bf Lemma}[section]
\newtheorem{Proposition}{\bf Proposition}[section]
\newtheorem{Corollary}{\bf Corollary}[section]
\newtheorem{Remark}{\bf Remark}[section]
\newtheorem{Example}{\bf Example}[section]
\newtheorem{Definition}{\bf Definition}[section]

\numberwithin{equation}{section}

\begin{document}

\title{
Existence of mild solutions for 
the Hamilton-Jacobi equation with critical fractional viscosity in the Besov spaces}
\author{
Tsukasa Iwabuchi\\
Department of Mathematics\\
Osaka City University and OCAMI\\
Sumiyoshi-ku, Osaka 558-8585, Japan
\\
\\
Tatsuki Kawakami\\
Department of Mathematical Sciences\\
Osaka Prefecture University\\
Sakai 599-8531, Japan}
\date{}
\maketitle
\begin{abstract}
We consider the Cauchy problem for the Hamilton-Jacobi equation 
with critical dissipation, 
$$
\partial_t u + (-\Delta)^{ 1/2} u = |\nabla u|^p, 
\quad x \in \mathbb R^N, t > 0, 
\qquad 
u(x,0) = u_0(x) , \quad x \in \mathbb R^N, 
$$
where $p > 1$ and $u_0 \in 
B^1_{r,1}(\mathbb R^N) \cap B^1_{\infty,1} (\mathbb R^N)$ 
with $r \in [1,\infty]$. 
We show that for sufficiently small $u_0 \in \dot B^1_{\infty,1}(\mathbb R^N)$, 
there exists a global-in-time mild solution. 
Furthermore, we prove that the solution behaves asymptotically like 
suitable multiplies of the Poisson kernel.  
\end{abstract}
\section{Introduction}
We consider the Hamilton-Jacobi equation with fractional viscosity, 
\begin{equation}
\label{eq:1.1}
\left\{
\begin{array}{ll}
\displaystyle{\partial_tu+(-\Delta)^{\alpha/2} u=|\nabla u|^p,} & x\in{\mathbb R}^N,\quad t>0,\vspace{5pt}\\
\displaystyle{u(x,0)=u_0(x)}, \quad 
& x\in{\mathbb R}^N,
\end{array}
\right.
\end{equation}
where $N\ge 1$, $\partial_t=\partial/\partial t$, $\nabla=(\partial_{x_1},\cdots,\partial_{x_N})$, $\partial_{x_j}=\partial/\partial x_{j}$ $(j=1,\cdots,N)$,
$\alpha\in(0,2]$, $p>1$ and $u_0$ is a nontrivial measurable function in ${\mathbb R}^N$. 
Here the operator $(-\Delta)^{\alpha/2}$, which called the L\'evy operator, is defined by the Fourier transform $\mathcal F$ such that
$$
(-\Delta)^{\alpha/2}f:={\mathcal F}^{-1}\big[|\xi|^\alpha{\mathcal F}[f]\big].
$$
In this paper we study the existence of global-in-time solutions to the problem \eqref{eq:1.1} with $\alpha=1$,
and investigate the asymptotics of solutions.
\vspace{5pt}

The problem \eqref{eq:1.1} with $\alpha=2$ is the well-known viscous Hamilton-Jacobi (VHJ) equation.
The VHJ equation possesses both mathematical and physical interest.
Indeed, in mathematical points of view,
it is the simplest example of a parabolic PDE with a nonlinearity depending only on the first order spatial derivatives of $u$,
and it describes a model for growing random interfaces,
which is known as the Kardar-Parisi-Zhang equation (see \cite{KPZ,KS}).
On the other hand, the problem \eqref{eq:1.1} with $\alpha\in(0,2)$ often appears in the context of  mathematical finance as Bellman equations
of optimal control of jump diffusion processes (see, for example, \cite{DI, I, JK1, JK2, S}).
\vspace{5pt}

The VHJ equation has been studied in many papers about various topics.
For the existence and uniqueness of solutions,
it is well known that, for any $u_0\in W^{1,\infty}({\mathbb R}^N)$,
the problem \eqref{eq:1.1} with $\alpha=2$ has a unique global-in-time mild solution, i.e., a solution of the integral equation
$$
u(t)=e^{t\Delta}u_0+\int_0^te^{(t-\tau)\Delta}|\nabla u(\tau)|^p\,d\tau,\qquad t>0,
$$ 
where $e^{t\Delta}$ denotes the convolution operator with the heat kernel
(see, for example, \cite{AB, BL, BSW, GGK}).
Furthermore, this solution is classical for positive time,
and by the maximum principle,
we see that, if $u_0\ge0$, then $u\ge0$, and if $u_0\le0$, then $u\le0$.
>From this property,
the nonlinearity $|\nabla u|^p$ behaves like a source term for nonnegative initial data and an absorption term for nonpositive initial data.
Similarly to the case 
of the semilinear heat equation $\partial_tu-\Delta u=\lambda|u|^{p-1}u$ with $\lambda=\pm1$,
the asymptotics of solutions to this equation is determined by the balance of effects  
from the diffusion term $\Delta u$ and the one from the nonlinearity $|\nabla u|^p$,
and there are many results on the asymptotic behavior of solutions.
See, for example, \cite{AB}--\cite{BSW}, \cite{GGK, KW, LS} and the references therein.
Among others, in \cite{BKL}, Benachour, Karch and Lanren\c{c}ot proved that,
for the case $u_0\in L^1({\mathbb R}^N)\cap W^{1,\infty}({\mathbb R}^N)$ with $u_0\not\equiv0$, 
the following hold.
\begin{itemize}
\item[(i)] Assume that $u_0\ge0$.
\begin{itemize}
\item[(a)] For the case $p\ge2$, there exists a limit
\begin{equation}
\label{eq:1.2}
C_*:=\lim_{t\to\infty}\int_{{\mathbb R}^N}u(x,t)\,dx=\int_{{\mathbb R}^N}u_0(x)\,dx+\int_0^\infty\int_{{\mathbb R}^N}|\nabla u(x,t)|^p\,dx\,dt
\end{equation}
such that
\begin{equation}
\label{eq:1.3}
\lim_{t\to\infty}t^{\frac{N}{2}(1-\frac{1}{q})+\frac{j}{2}}\|\nabla^j[u(t)-C_*G(t)]\|_{L^q({\mathbb R}^N)}=0,\qquad q\in[1,\infty],\quad j=0,1,
\end{equation}
where $G(x,t)$ is the heat kernel.

\item[(b)] For the case $p\in(p_c,2)$ with $p_c:=(N+2)/(N+1)$,
there exists a positive constant $\varepsilon=\varepsilon(N,p)$ such that, if
$$
\|u_0\|_{L^1({\mathbb R}^N)}\|\nabla u_0\|_{L^\infty({\mathbb R}^N)}^{(N+1)p-(N+2)}<\varepsilon,
$$
then \eqref{eq:1.3} holds ture.
\end{itemize}
\item[(ii)] Assume that $u_0\le0$.
For any $p>p_c$, \eqref{eq:1.3} holds true.
\end{itemize}
For the case $\alpha\in(1,2)$,
Karch and Woyczy\'nski \cite{KW} studied similar topics.
They showed  that, for any $u_0\in W^{1,\infty}({\mathbb R}^N)$, 
the problem \eqref{eq:1.1} with $\alpha\in(1,2)$ has a unique global-in-time mild solution.
Furthermore, for the case $p>(N+\alpha)/(N+1)$,
they proved that there exists a (mild) solution which behaves asymptotically like suitable multiples of the kernel of the integral equation.
For notions of another weak solutions,
Droniou and Imbert \cite{DI} constructed a unique global-in-time viscosity solution in $W^{1,\infty}({\mathbb R}^N)$ for the case $\alpha\in(0,2)$
(see also \cite{I,Si}).

On the other hand,
the case $\alpha=1$ is completely different from the case $\alpha\in(1,2]$. 
In fact, for the case $\alpha\in(0,2]$, the semigroup $e^{-t(-\Delta)^{\alpha/2}}$ satisfies the following decay estimates
$$
\|\nabla^je^{-t(-\Delta)^{\alpha/2}} f\|_{L^q({\mathbb R}^N)}\le Ct^{-\frac{N}{\alpha}(1-\frac{1}{q})-\frac{j}{\alpha}}\|f\|_{L^1({\mathbb R}^N)},
\qquad q\in[1,\infty],\quad j=0,1,
$$
for all $t>0$ (see, for example, \cite{IKK2}).
For the case $\alpha\in(1,2]$, since $t^{-1/\alpha}$ is integrable locally,
we can easily prove the existence of local-in-time mild solutions in $W^{1,\infty}({\mathbb R}^N)$
(see \cite[Proposition~3.1]{KW}).
However, for the case $\alpha=1$, since $t^{-1}$ is not integrable, we need to impose the regularity of one order derivative on the solution.
In this sense the value $\alpha=1$ is critical.
Similar situation appears in the fractional Burgers equation, 
\begin{equation}
\label{eq:1.4}
\partial_t u+u\partial_xu+(-\partial_{xx})^{\alpha/2}u=0,\qquad x\in\mathbb R,\quad t>0.
\end{equation}
For \eqref{eq:1.4}, the value $\alpha=1$ is a threshold for the occurrence of singularity in finite time or the global regularity
(see \cite{AIK, DDL, DGV, KNS}).
In \cite{Iw}, the first author of this paper studied \eqref{eq:1.4} with $\alpha=1$, 
and constructed a small global-in-time mild solution in the Besov space $\dot{B}^0_{\infty,1}(\mathbb R)$ 
which is the critical space under the scaling invariance (see also \cite{MW}).
Furthermore, he proved that, for small initial data in $L^1(\mathbb R)\cap \dot{B}^0_{\infty,1}(\mathbb R)$,
the corresponding solution behaves like the Poisson kernel as $t\to\infty$.
\vspace{5pt}

In this paper, modifying the argument in \cite{Iw},
we show that there exists a global-in-time mild solution of the problem \eqref{eq:1.1} with $\alpha=1$ in the critical Besov space.
Furthermore, we prove that global-in-time solutions 
with some suitable decay estimates 
behave asymptotically like suitable multiples of the Poisson kernel.
\vspace{5pt}

We introduce some notations.
Throughout this paper we put $\mathcal {L}:=-(-\Delta)^{1/2}$ for simplicity.
Let $P_t$ be the Poisson kernel, that is,
$$
P_t(x):=t^{-N}P(x/t),\qquad x\in{\mathbb R}^N,\quad t>0,
$$
where
$P$ is defined by
$$
P(x):=\mathcal{F}^{-1}\big[e^{-|\xi|}\big](x)=c_N(1+|x|^2)^{-(N+1)/2},\qquad x\in{\mathbb R}^N,
$$
and $c_N$ is a constant chosen so that 
\begin{equation}
\label{eq:1.5}
\int_{{\mathbb R}^N}P(x)\,dx=1. 
\end{equation}
Then, for all $t>0$, $e^{t\mathcal{L}}$ denotes the convolution operator with $P_t$, that is,
\begin{equation}
\label{eq:1.6}
[e^{t\mathcal{L}}f](x):=\int_{{\mathbb R}^N}P_t(x-y)f(y)\,dy,\qquad x\in{\mathbb R}^N,\quad t>0,
\end{equation}
and $f$ is a measurable function.
For $q \in [1,\infty]$, we denote by $\|\cdot\|_{L^q}$ the usual norm of $L^q:= L^q({\mathbb R}^N)$.
Furthermore, for $s\in{\mathbb R}$, $q\in [1,\infty]$ 
and $\sigma \in (0,\infty]$, 
we denote by $\|\cdot\|_{B^s_{q,\sigma}}$ and $\|\cdot\|_{\dot B^s_{q,\sigma}}$ 
the usual norm of inhomogeneous and homogeneous Besov spaces 
$B^s_{ q,\sigma}:=B^s_{q,\sigma}({\mathbb R}^N)$ and $\dot B^s_{q,\sigma}:=\dot B^s_{q,\sigma}({\mathbb R}^N)$, respectively.
(See Section 2 for more precise details.)
\vspace{5pt}

Now we are ready to state the main result of this paper.
We consider the integral equation corresponding to \eqref{eq:1.1} with $\alpha=1$, that is,
\begin{equation}
\label{eq:1.7}
u(t) = e^{t{\mathcal L}}u_0+\int_0^t e^{(t-\tau)\mathcal L}|\nabla u(\tau)|^p\,d\tau,\qquad t\ge0,
\end{equation} 
and obtain the following result. 
%
\begin{Theorem}
\label{Theorem:1.1}
Let $N \geq 1$, $p > 1$ and $r\in[1,\infty]$.
Assume $u_0 \in B^1_{q,1}$ for all $q\in[r,\infty]$.
Then the following hold.  
\begin{itemize}
\item[\rm(i)]
There exists a positive constant $\delta=\delta(N,p)$ such that, 
if
\begin{equation}
\label{eq:1.8}
\| u_0 \|_{\dot B^1_{\infty,1}} \le \delta,
\end{equation}
then there exists a unique global-in-time solution $u$ of \eqref{eq:1.7} satisfying
$$
u\in C([0,\infty), B^1_{q,1})\cap 
L^1 (0,\infty ; \dot B^2_{q,1}),
$$
\begin{equation}
\label{eq:1.9}
  \sup_{t\ge0}(1+t)^{N(\frac{1}{r}-\frac{1}{q})+j}\| \nabla^ju(t) \|_{L^q}<\infty,
\end{equation}
\begin{equation}
\label{eq:1.10}
\int_0^\infty t^{N(\frac{1}{r}-\frac{1}{q})+\frac{1}{p}} \| u(t) \|_{\dot B^{2}_{q,1}}\,dt<\infty,
\end{equation}
for all $q\in[r,\infty]$ and $j=0,1$. 
\vspace{5pt}
\item[\rm(ii)]
Let $v$ be a global-in-time solution of \eqref{eq:1.7} satisfying \eqref{eq:1.9} and \eqref{eq:1.10}.
Then for any $j\in\{0,1\}$, the following hold.  
\begin{itemize}
  \item[{\rm (a)}] If $1<r<\infty$, then 
  \begin{equation}
  \label{eq:1.11}
  t^{N(\frac{1}{r}-\frac{1}{q})+j}\left\|\nabla^j[v(t)-e^{t{\mathcal L}}u_0]\right\|_{L^q}
  =
  \left\{
  \begin{array}{ll}
  O(t^{-\frac{N}{r}(r-1)})
  &
  \mbox{if}\quad p\ge r,
  \vspace{5pt}\\
  O(t^{-\frac{N}{r}(p-1)})
  &
  \mbox{if}\quad p<r,
  \end{array}
  \right.
  \end{equation}
  as $t\to\infty$, for any $q\in[r,\infty]$.
  \item[{\rm (b)}] If $r=1$, then the limit $C_*$ given in \eqref{eq:1.2} exists and 
  \begin{equation}
  \label{eq:1.12}
  \lim_{t\to\infty}t^{N(1-\frac{1}{q})+j}\left\|\nabla^j[v(t)-C_*P_{t+1}]\right\|_{L^q}=0,\qquad 1\le q\le\infty.
  \end{equation}
\end{itemize}
\end{itemize}
\end{Theorem}
\begin{Remark}
{\rm(i)}
Let $u$ be a mild solution $u$ of \eqref{eq:1.1} with $\alpha=1$, i.e., solution of \eqref{eq:1.7}. 
For any $\lambda>0$, put 
\begin{equation}
\label{eq:1.13}
u_\lambda(x,t):=\lambda^{-1} u(\lambda x,\lambda t),\qquad
u_{0,\lambda}(x):=\lambda^{-1}u_0(\lambda x).
\end{equation}
Then the function $u_\lambda$  is also a solution of \eqref{eq:1.1} with $\alpha=1$ and the initial function $u_{0,\lambda}$ satisfying
\begin{equation}
\label{eq:1.14}
C^{-1}\|u_0\|_{\dot B^1_{\infty,1}} 
\leq \|u_{0,\lambda}\|_{\dot B^1_{\infty,1}} 
\leq C \|u_0\|_{\dot B^1_{\infty,1}}, 
\end{equation}
where $C $ is a positive constant independent of $\lambda$. 
This means that the condition~\eqref{eq:1.8} is invariant with respect to the similarity transformation~\eqref{eq:1.13}.
This is the reason why we say that $\dot B^1_{\infty,1}$ is the critical Besov space with respect to \eqref{eq:1.1}. 
\vspace{5pt}
\newline
{\rm(ii)}
In the assertion~{\rm(ii)} of Theorem~$\ref{Theorem:1.1}$, if we only consider the case $j=0$,
then we can remove the assumption that the solution $u$ satisfies \eqref{eq:1.10}.
See Section~$5$.
\vspace{5pt}
\newline
{\rm(iii)} 
As is seen from our proof, it is possible to replace 
\eqref{eq:1.10} with 
$$
\int_0^\infty t^{N(\frac{1}{r} - \frac{1}{q}) + \beta} 
  \| u(t) \|_{\dot B^2_{q,1}} \, dt 
  < \infty, 
$$
where $1/p < \beta < 1$. 
We also note focusing on the linear part 
that for $\beta = 1$, the maximal regularity estimate 
and the embedding implies that 
$$
\int_0^\infty t^{N(\frac{1}{r} - \frac{1}{q}) + 1} 
  \| e^{t\mathcal L} u_0 \|_{\dot B^2_{q,1}} \, dt 
  \leq C \| u_0 \|_{\dot B^0_{r,1}} ,
  \quad 
  B^1_{r,1} \not \subset \dot B^0_{r,1}, 
$$
and one can not expect the time decay with $\beta = 1$ for initia data 
in $B^1_{r,1}$. 
Therefore, the expected maximal decay order is given as the case $\beta < 1$ 
expect for $\beta = 1$, 
and 
the case $\beta = 1/p$ is a sufficient decay to prove the asymptotic behavior. 
\vspace{5pt}
\newline
{\rm(iv)} 
By the embedding 
$B^1_{\infty,1} \hookrightarrow C^1$ 
and $(-\Delta)^{1/2} f \in C(\mathbb R^N)$ for 
$ f \in B^1_{\infty,1}$, 
the solution $u$ in Theorem~{\rm\ref{Theorem:1.1}}-{\rm (i)}
satisfies the problem \eqref{eq:1.1} in the classical sense. 
We also see that $u(t)$ is in the class $C^2$ for almost every $t$ 
since $u \in L^1(0,\infty ; \dot B^2_{\infty,1})$. 
Compared with the results~{\rm\cite{DI,I,Si}}, 
our framework in the Besov spaces is the one with higher regularity 
than theirs, 
since their initial data are in $W^{1,\infty}$ and solutions 
are considered in the sense of viscosity solutions 
and $B^1_{\infty,1}\hookrightarrow W^{1,\infty}$. 
\vspace{5pt}
\newline
{\rm(v)} 
In Theorem~{\rm\ref{Theorem:1.1}-(ii)-(a)}, 
it is possible to prove that 
$$
  \lim_{t\to\infty}t^{N(\frac{1}{r}-\frac{1}{q})+j}\left\|\nabla^ju(t)\right\|_{L^q}=0
$$
since one can show 
$t^{N(\frac{1}{r}-\frac{1}{q})+j}\left\|\nabla^j e^{t\mathcal L}u_0 \right\|_{L^q} 
 =o(1)$ as $t\to\infty$ for any $u_0 \in L^r$ by the density argument due to $C_0^\infty \subset L^r$. 
\end{Remark}
\vspace{5pt}

This paper is organized as follows.
In Section 2, 
we give the definition of the Besov spaces, its properties 
and estimates for the nonlinearity $|\nabla u|^p$. 
We also introduce the linear estimates for $e^{t\mathcal L} f$ 
in the Lebesgue spaces and the Besov spaces. 
Sections 3 and 4 
are devoted to the proof of the assertions 
{\rm (i)} and {\rm (ii)} in Theorem \ref{Theorem:1.1}, 
respectively. 

\section{Preliminary}
In this section 
we prove some estimates in the Besov spaces
and  recall some preliminary results on $e^{t{\mathcal L}}f$. 
In what follows, 
for any two nonnegative functions 
$f_1$ and $f_2$ on a subset $D$ of $[0,\infty)$, 
we say 
$$
f_1(t)\preceq f_2(t),\qquad t\in D
$$ 
if there exists a positive constant $C$ such that 
$f_1(t)\le Cf_2(t)$ for all $t\in D$. 
In addition, we say  
$$
f_1(t)\asymp f_2(t),\qquad t\in D
$$ 
if $f_1(t)\preceq  f_2(t)$ and $f_2(t)\preceq f_1(t)$ for all $t\in D$. 
We denote the function spaces of rapidly decreasing functions by  $\mathcal{S}(\mathbb R^N)$ 
and tempered distributions by $\mathcal{S}'(\mathbb R^N)$. 
We define $\mathcal Z (\mathbb R^N)$ by 
$$
\mathcal Z(\mathbb R^N) 
:= 
\Big\{ f \in \mathcal S (\mathbb R^N) 
\, \Big| \, 
\int_{\mathbb R^N} x^\alpha f(x) dx = 0  \,\, 
\text{ for all } \alpha \in (\{ 0 \} \cup \mathbb N) ^N 
\Big\} 
$$
with the topology of $\mathcal S (\mathbb R^N)$,
and $\mathcal Z'(\mathbb R ^N)$ by the topological dual of $\mathcal Z(\mathbb R^N)$. 
We first give the definition of the inhomogeneous and homogeneous Besov spaces 
(see Triebel~\cite{Tri_1983}).
\vspace{5pt}

\begin{Definition}
\label{Definition:2.1}
Let $\phi \in \mathcal S (\mathbb R^N) $ satisfy 
$$
{\rm supp \,} \mathcal F [\phi]\subset \{ \, \xi \in \mathbb R^N \, | \, 
   2^{-1} \leq |\xi| \leq 2 \}, 
\quad 
      \sum _{j \in \mathbb Z} \mathcal F [\phi] (2^{-j} \xi) = 1 
   \,\,\, \text{for any } \xi \in \mathbb R^N \setminus  \{ 0 \},  
$$
Let $\{ \phi_j \}_{j \in \mathbb Z}$ and $\psi$ be defined by 
\begin{equation}\notag 
\phi_j (x) := 2^{Nj} \phi (2^j x), 
\quad 
\psi (x)
= \mathcal F^{-1}
\Big[ 1 - \sum _{j \geq 1} \mathcal F[\phi_j ]
\Big](x)
\end{equation}
For $s \in \mathbb R$, $q \in [1,\infty]$ and $\sigma \in (0,\infty]$, 
we define the following. 

\begin{itemize}
\item[{\rm (i)}] 
The inhomogeneous Besov space $B^s_{q,\sigma}$ is defined by 
\begin{equation}\notag 
  B^{s}_{q,\sigma} := 
 \big\{ \,\, 
  u \in \mathcal S' (\mathbb R ^N) 
  \,\, \big| \,\, 
  \| u \|_{B^{s}_{q,\sigma}} 
   < \infty
 \,\, \big\} , 
\end{equation}
where 
\begin{equation}\notag
   \| u \|_{B^{s}_{q,\sigma}} := 
\begin{cases} 
 \displaystyle 
   \| \psi * u \|_{L^q} + 
     \Big\{
      \sum _{j \geq 1} 
       \big( 
        2^{js} \| \phi_j * u \|_{L^q }  
       \big)^\sigma 
     \Big\} ^{1/\sigma} 
 & \text{if} \quad 0< \sigma < \infty, 
 \\ 
 \displaystyle 
   \| \psi * u \|_{L^q} + 
  \sup _{j \geq 1}  
        2^{js} \| \phi_j *  u \|_{L^q }  
 & \text{if} \quad \sigma = \infty. 
\end{cases}
\end{equation}

\item[{\rm (ii)}] 
The homogeneous Besov space $\dot B^s_{q,\sigma }$ is defined by 
\begin{equation}\notag 
 \dot B^{s}_{q,\sigma} := 
 \big\{ \,\, 
  u \in \mathcal Z' (\mathbb R^N)
  \,\, \big| \,\, 
  \| u \|_{\dot B^{s}_{q,\sigma}} 
   < \infty
 \,\, \big\} , 
\end{equation}
where 
\begin{equation}\notag 
   \| u \|_{\dot B^{s}_{q,\sigma}} := 
\begin{cases} 
 \displaystyle 
     \Big\{
      \sum _{j \in \mathbb Z} 
       \big( 
        2^{js} \| \phi_j * u \|_{L^q }  
       \big)^\sigma 
     \Big\} ^{1/\sigma} 
 & \text{if} \quad 0 < \sigma < \infty, 
 \\ 
 \displaystyle 
  \sup _{j \in \mathbb Z}  
        2^{js} \| \phi_j *  u \|_{L^q }  
 & \text{if} \quad \sigma = \infty.    
\end{cases}
\end{equation}
\end{itemize}
\end{Definition}
\begin{Remark} 
\label{Remark:2.1}
It is known that 
$\mathcal Z(\mathbb R^N) \subset \mathcal S(\mathbb R^N) 
\subset \mathcal S'(\mathbb R^N) \subset \mathcal Z'(\mathbb R^N)$ 
and 
$\mathcal Z'(\mathbb R^N ) \simeq \mathcal S' (\mathbb R^N)/ \mathcal P (\mathbb R^N)$, 
where $\mathcal P (\mathbb R^N)$ is the set of all polynomials, 
and the homogeneous Besov spaces can also be considered as subspaces of the quotient space 
$\mathcal S'(\mathbb R^N) / \mathcal P (\mathbb R^N)$. 
Then we use the following equivalence, 
which is due to the argument by e.g. Kozono and Yamazaki~{\rm\cite{KoYa-1994}}, 
for the nonlinear term in \eqref{eq:1.1} to construct solutions in the homogeneous spaces 
with $u(t) \in \mathcal S'(\mathbb R^N)$.
If $s<n/q$ or $(s,\sigma) = (n/q,1)$, then the homogeneous Besov space $\dot B^s_{q,\sigma}$
is regarded as 
\begin{equation}\notag 
\Big\{ \, 
 u \in \mathcal S' (\mathbb R^N) 
 \, \Big| \, 
 \| u \|_{\dot B^s_{q,\sigma}} < \infty , \,
 u= \sum _{j \in \mathbb Z} \phi_j * u \text{ in } 
 \mathcal S'(\mathbb R^N)
\,\Big\}. 
\end{equation}
Hence, $u \in \dot B^s_{q,\sigma}$ can be regarded as an element of 
$\mathcal S '(\mathbb R^N)$. 
We also see from the analogous argument to theirs that 
$\nabla u $ can be regarded as an element of $\mathcal S'(\mathbb R^N)$ 
if $u \in \dot B^s_{\infty,1}$ with $s \leq 1$. 
This is used for the nonlinear term $|\nabla u|^p$ 
when we construct global solutions. 
\end{Remark}

\noindent
Next we give some interpolation inequalities in the Besov spaces.
\begin{Lemma}
\label{Lemma:2.1} 
Let $s \in \mathbb R$, $\alpha , \beta  > 0$, 
$q \in [1,\infty]$ and $\sigma \in (0,\infty]$. 
Then it holds that 
\begin{equation}\label{eq:2.1}
\| f \|_{\dot B^s_{q,\sigma}} 
\preceq 
\| f \|_{\dot B^{s+\alpha}_{q,\infty}} ^{\frac{\beta}{\alpha + \beta}} 
       \| f \|_{\dot B^{s-\beta}_{q,\infty}} ^{\frac{\alpha}{\alpha + \beta}}
\end{equation}
for all $f \in \dot B^{s+\alpha}_{q,\infty}
 \cap \dot B^{s-\beta}_{q,\infty} $. 
\end{Lemma}
{\bf Proof.}
The estimate \eqref{eq:2.1} is known for the case $1 \leq \sigma \leq \infty$ 
by the result of Machihara-Ozawa~{\rm\cite{MO-2003}}.  
The case $0 < \sigma < 1$ follows from the analogous argument to their proof, 
thus the proof is left to readers. 
$\Box$\vspace{7pt}

\noindent
The following proposition is on the equivalence between the norm of 
the Besov spaces defined as above and that by differences (see Triebel~{\rm\cite{Tri_1983}}). 
\begin{Proposition}
\label{Proposition:2.1}
Let $s>0$, $q \in [1,\infty]$ and $\sigma \in (0,\infty ]$. 
If $M \in \mathbb N $ satisfies $M > s$, 
then there holds that 
\begin{equation}\label{eq:2.2}
\| f \|_{\dot B^s_{q,\sigma}}  
\asymp
\Big\{ 
\int_{\mathbb R^N}
 \Big( 
  |\eta|^{-s}
  \sup _{|y| \leq |\eta|} \| \triangle_y ^M f \|_{L^q}
 \Big) ^{\sigma}
  \frac{d\eta}{|\eta|^ N}
\Big\} ^{ 1/\sigma} 
\end{equation}
for all $f \in \dot B^s_{q,\sigma}$, 
where $\triangle_y f (x) := f(x+y) - f(x)$ 
and $\triangle _y ^M f := ( \triangle _y )^M f$. 
\end{Proposition}

\noindent
By using Proposition~\ref{Proposition:2.1} we have the following.
\begin{Lemma}
\label{Lemma:2.2} 
Let $p$, $s$ and $\varepsilon$ satisfy 
$p > 1$, $0 < s < \min \{ 2,p \}$ and $ 0 < \varepsilon < \min\{ 1 , p-1 \} $. 
Then, for any $q\in[1,\infty]$,
\begin{eqnarray}\label{eq:2.3}
&&\| |f|^p \|_{\dot B^s_{q,1}} 
 \preceq \| f \|_{\dot B^0_{\infty,1}}^{p-1} 
       \| f \|_{\dot B^s_{q,1}} , 
\\ 
\label{eq:2.4}
&&
\begin{aligned}
&
\| |f|^p - |g|^p \|_{\dot B^\varepsilon_{\infty,1}} 
\\
&
\preceq \big( \| f \|_{\dot B^0_{\infty,1}} ^{p-1} + \| g \|_{\dot B^0_{\infty,1}}^{p-1} 
       \big) 
       \| f-g \|_{\dot B^\varepsilon_{q,1}}
\\
&\quad
+
\left\{
\begin{array}{l}
\displaystyle{
\Big(
             \| f \|_{\dot B^0_{\infty, 1}} ^{p-1-\varepsilon}
             \| f \|_{\dot B^1_{\infty, 1}} ^{\varepsilon}
            +
             \| g \|_{\dot B^0_{\infty, 1}} ^{p-1-\varepsilon}
             \| g \|_{\dot B^1_{\infty, 1}} ^{\varepsilon}
       \Big)
\| f-g \|_{\dot B^0_{q,1}}},
\\
\hspace{8.7cm} 
\text{if}\quad 1 < p < 2, 
\vspace{5pt}\\
\displaystyle{
\Big( \| f \|_{\dot B^0_{\infty, 1}} ^{p-2} + \| g \|_{\dot B^0_{\infty, 1}} ^{p-2} 
       \Big)   
       \Big( \| f \|_{\dot B^0_{\infty,1}} ^{1-\varepsilon} 
             \| f \|_{\dot B^1_{\infty,1}} ^\varepsilon
           + \| g \|_{\dot B^0_{\infty,1}} ^{1-\varepsilon} 
             \| g \|_{\dot B^1_{\infty,1}} ^\varepsilon
       \Big) 
       \| f-g \|_{\dot B^0_{q,1}}},
       \\
       \hspace{8.7cm}
\text{if}\quad p \geq 2 , 
\end{array}
\right.
\end{aligned}
\end{eqnarray}
for all 
$f,g \in \dot B^0_{\infty,1}  \cap \dot B^{1}_{\infty,1} 
 \cap \dot B^0_{q,1} \cap \dot B^{\max\{s,\varepsilon\}}_{q,1}$. 
\end{Lemma}

\noindent
In order to prove this lemma,
we use the following fundamental inequality.
\begin{Lemma}
\label{Lemma:2.3}
Let $p>1$.
Then, for any $A, B, C, D\in{\mathbb R}$,
\begin{equation}
\label{eq:2.5}
\begin{split}
& 
\big||A|^p- |B|^p - (|C|^p-|D|^p) \big| 
\\
& 
\preceq 
\big( |
|C|^{p-1} + |D|^{p-1} \big) 
|A-B - (C-D)| 
\\
&
+ 
\left\{
\begin{array}{ll}
\displaystyle{
\big( |A-C|^{p-1} + |B-D|^{p-1} \big) |A-B|}, 
&\hspace{-21pt} \text{if}\,\,\, 1 < p < 2, 
\vspace{5pt}\\
\displaystyle{
\big( |A|^{p-2} + |B|^{p-2} + |C|^{p-2} + |D|^{p-2} \big) 
\big( |A-C| + |C-D| \big) |A-B|},
& \text{if}\,\,\, p \geq 2.
\end{array}
\right.
\end{split}
\end{equation}
\end{Lemma}
{\bf Proof.}
Let $p>1$ and $A, B, C, D\in{\mathbb R}$.
It follows from the fundamental theorem of calculus that 
\begin{equation}
\label{eq:2.6}
\begin{split}
& |A|^p - |B|^p 
      - \big( |C|^p - |D|^p 
        \big) 
\\
& 
 = \int_0^1 
  \Big\{ 
      \partial_\theta 
      \big| \theta A + (1-\theta) B 
      \big|^p 
      - 
      \partial_\theta 
      \big| \theta C + (1-\theta) D 
      \big|^p
  \Big\} 
  d\theta 
\\
& 
 = \int_0^1 
  \Big\{ 
      \big| \theta A + (1-\theta) B
      \big|^{p-2} 
      \big( \theta A + (1-\theta) B
      \big) 
      \big( A-B 
      \big)
\\
& 
\qquad \qquad 
      - 
      \big| \theta C + (1-\theta) D 
      \big|^{p-2}
      \big( \theta C + (1-\theta) D 
      \big)  
      \big(C-D
      \big)
  \Big\} 
  d\theta 
\\
& 
 = \int_0^1 
  \Big[
  \Big\{
      \big| \theta A + (1-\theta) B 
      \big|^{p-2} 
      \big( \theta A + (1-\theta) B
      \big) 
\\
& 
\qquad \qquad 
- 
      \big| \theta C + (1-\theta) D 
      \big|^{p-2} 
      \big( \theta C + (1-\theta) D
      \big) 
  \Big\}
  \big( A - B 
      \big)
\\
& 
\qquad \quad 
      + 
      \big| \theta C + (1-\theta) D 
      \big|^{p-2}
      \big( \theta C + (1-\theta) D 
      \big) 
      \Big( A - B - \big( C - D \big) 
      \Big)
  \Big]
  d\theta.
\end{split}
\end{equation}
Furthermore, we have
\begin{equation*} 
\big| |E|^{p-2} E - |F|^{p-2} F \big| 
\leq C
\begin{cases}
|E-F|^{p-1} ,
& \text{if } 1 < p < 2 ,
\\
(|E|^{p-2} + |F|^{p-2}) |E-F| , 
& \text{if } p \geq 2,
\end{cases}
\end{equation*}
for any $E, F \in \mathbb R$.
This together with \eqref{eq:2.6} yields \eqref{eq:2.5}.
Thus Lemma~\ref{Lemma:2.3} follows. 
$\Box$\vspace{7pt}
\noindent
{\bf Proof of Lemma~\ref{Lemma:2.2}.}
For the proof of \eqref{eq:2.3}, 
we utilize the equivalent norm \eqref{eq:2.2} of the Besov spaces 
$\dot B^s_{q,1} $ by differences, and 
it suffices to estimate the following 
$$
\int_{\mathbb R^N}
\Big( |\eta|^{-s} 
\sup _{ |y| \leq |\eta|} 
 \big\| \triangle _y^2 |f|^p 
 \big\|_{L^q} 
\Big) 
\frac{d\eta}{|\eta|^N} . 
$$
In order to estimate $\triangle ^2 _y |f|^p$, 
we apply Lemma~\ref{Lemma:2.2}.
Put 
$$
A =f(x+2y),\qquad B = C = f(x+y),\qquad D = f(x),
$$ 
and we note that 
$$
A-B = (\triangle_y f) (x+y),\qquad 
C-D = (\triangle_y f) (x), \qquad
A-B-(C-D) = (\triangle_y^2 f) (x).
$$ 
In the case $ 1 < p < 2$, 
by \eqref{eq:2.5} and the H\"older inequality we have
\begin{equation}\label{eq:2.7}
\begin{split}
& \| \triangle^2_y |f|^p \|_{L^q} 
\\ & 
\preceq 
 \| f \|_{L^\infty} ^{p-1} \| \triangle ^2_y f \|_{L^q} 
 + \big( \| | \triangle _y f (\cdot + y)|^{p-1} \|_{L^{\frac{pq}{p-1}}}  
       + \| | \triangle _y f |^{p-1} \|_{L^{\frac{pq}{p-1}}}  
   \big)
   \| \triangle _y f \|_{L^{pq}} 
\\
& 
\preceq 
 \| f \|_{L^\infty} ^{p-1} \| \triangle ^2_y f \|_{L^q} 
 +  \| \triangle _y f \|_{L^{pq}} ^p 
\\
& 
=: I_1(y) + I_2(y) . 
\end{split}
\end{equation}
On the estimate of $I_1$, we get 
\begin{equation}\label{eq:2.8}
\begin{split}
&
\int_{\mathbb R^N} 
|\eta|^{-s} 
\sup _{ |y| \leq |\eta|} I_1(y) \frac{d\eta}{|\eta|^N}
\\
& 
\preceq \| f \|_{L^\infty}^{p-1}
\int_{\mathbb R^N} 
\Big( |\eta|^{-s} 
\sup _{ |y| \leq |\eta|} \| \triangle _y^2 f \|_{L^q}
\Big)
\frac{d\eta}{|\eta|^N}
\preceq \| f \|_{\dot B^0_{\infty,1}} ^{p-1} 
     \| f \|_{\dot B^s_{q,1}}. 
\end{split}
\end{equation}
On the estimate of $I_2$,
applying the H\"older inequality and the embedding $\dot B^0_{\infty,1}\hookrightarrow \dot B^0_{\infty,\infty}$,
we have
\begin{equation}\label{eq:2.9}
\begin{split}
\int_{\mathbb R^N} |\eta|^{-s} 
\sup _{ |y| \leq |\eta|} I_2(y) \frac{d\eta}{|\eta|^N}
& 
= 
\int_{\mathbb R^N} 
\Big( |\eta|^{-\frac{s}{p}} \sup _{ |y| \leq |\eta|} \| \triangle _y f \|_{L^{pq}} 
\Big)^p
 \frac{d\eta}{|\eta|^N}
\\
& 
\preceq \| f \|_{\dot B^{\frac{s}{p}}_{pq,p}} ^p
\preceq \Big( \| f \|_{\dot B^0_{\infty,\infty}} ^{1-\frac{1}{p}} 
              \| f \|_{\dot B^s_{q,1}}^{\frac{1}{p}}
        \Big) ^{p}
\preceq \| f \|_{\dot B^0_{\infty,1}} ^{p-1} 
     \| f \|_{\dot B^s_{q,1}}. 
\end{split}
\end{equation}
In the case $p \geq 2$, 
by \eqref{eq:2.5} and the H\"older inequality again we obtain
\begin{equation}
\label{eq:2.10}
\begin{split}
\| \triangle^2_y f \|_{L^q} 
& 
\preceq 
 \| f \|_{L^\infty} ^{p-1} \| \triangle ^2_y f \|_{L^q} 
 + \| f \|_{L^\infty} ^{p-2}
   \big( \| \triangle _y f (\cdot + y) \|_{L^{2q}} + \| \triangle _y f  \|_{L^{2q}} \big) 
\\
&  
\quad \times \big( \| \triangle _y f (\cdot + y) \|_{L^{2q}} + \| \triangle _y f  \|_{L^{2q}} \big) 
\\
& 
\preceq 
 \| f \|_{L^\infty} ^{p-1} \| \triangle ^2_y f \|_{L^q} 
 +  \| f  \|_{L^\infty} ^{p-2} \| \triangle _y f \|_{L^{2q}}^2 
\\
& 
=: I_1(y) + I_3(y) . 
\end{split}
\end{equation}
For the estimate of $I_3$, 
it follows from the same estimate as \eqref{eq:2.9} with taking $p=2$ 
for $\| \triangle _y f \|_{L^{2q}}$ that 
\begin{equation}\label{eq:2.11}
\begin{split}
\int_{\mathbb R^N} 
|\eta|^{-s} 
\sup _{ |y| \leq |\eta|} I_3(y) \frac{d\eta}{|\eta|^N}
& 
\preceq 
\| f \|_{\dot B^0_{\infty,1}} ^{p-2}\Big( 
\| f \|_{\dot B^0_{\infty,1}}     \| f \|_{\dot B^s_{q,1}}\Big)
= 
\| f \|_{\dot B^0_{\infty,1}} ^{p-1} 
     \| f \|_{\dot B^s_{q,1}}. 
\end{split}
\end{equation}
Combining \eqref{eq:2.7}, \eqref{eq:2.8}, \eqref{eq:2.9}, \eqref{eq:2.10} and \eqref{eq:2.11}, 
the estimate \eqref{eq:2.3} holds.

For the proof of \eqref{eq:2.4}, 
we also utilize the equivalent norm \eqref{eq:2.2} of the Besov space 
$\dot B^\varepsilon_{\infty,1} $ by differences, and 
it suffices to estimate the following  
$$
\int_{\mathbb R^N} 
\Big( |\eta|^{-\varepsilon} 
\sup _{ |y| \leq |\eta|} 
 \big\| \triangle _y ( |f|^p - |g|^p) 
 \big\|_{L^q} 
\Big) 
\frac{d\eta}{|\eta|^N} . 
$$
In order to estimate $\triangle _y ( |f|^p - |g|^p )$, 
we also apply Lemma~\ref{Lemma:2.2}.
Put 
$$
A =f(x+y),\qquad 
B = g(x+y),\qquad 
C = f(x),\qquad
D = g(x),
$$ 
we note that 
$$
A-C = (\triangle_y f) (x),\qquad 
B-D = (\triangle_y g) (x),\qquad 
A-B-(C-D) = \triangle_y (f-g).
$$ 
In the case $ 1 < p < 2$, 
by \eqref{eq:2.5} and the H\"older inequality we have
\begin{equation}
\label{eq:2.12}
\begin{split}
\| \triangle_y (|f|^p - |g|^p ) \|_{L^q} 
& 
\preceq 
 \big( \| f \|_{L^\infty} ^{p-1} + \| g \|_{L^\infty} ^{p-1} \big) 
       \| \triangle _y (f-g) \|_{L^q} 
\\
& \quad 
 +  \big( \| \triangle _y f \|_{L^\infty}^{p-1} + \| \triangle _y g  \|_{L^\infty}^{p-1} \big) 
   \| f-g \|_{L^q} 
\\
& 
=: J_1(y) + J_2(y) . 
\end{split}
\end{equation}
On the estimate of $J_1$, we get 
\begin{equation}
\label{eq:2.13}
\begin{split}
\int_{\mathbb R^N} 
|\eta|^{-\varepsilon} 
\sup _{ |y| \leq |\eta|} J_1(y) \frac{d\eta}{|\eta|^N}
& 
\preceq \big( \| f \|_{L^\infty}^{p-1} + \| g \|_{L^\infty}^{p-1} \big)
\int_{\mathbb R^N} 
\Big( 
|\eta|^{-\varepsilon} 
\sup _{ |y| \leq |\eta|} \| \triangle _y (f-g) \|_{L^q}
\Big) 
\frac{d\eta}{|\eta|^N}
\\
& 
\preceq \big( \| f \|_{\dot B^0_{\infty,1}}^{p-1} + \| g \|_{\dot B^0_{\infty,1}}^{p-1} \big)
     \| f-g \|_{\dot B^\varepsilon_{q,1}}. 
\end{split}
\end{equation}
On the estimate of $J_2$, 
by \eqref{eq:2.1} and the embeddings 
$\dot B^{s}_{q,1} 
 \hookrightarrow \dot B^{s}_{q,\infty} $ ($s = 0,1$) 
 and $\dot B^0_{q,1} \hookrightarrow L^q$
we have
\begin{equation}
\label{eq:2.14}
\begin{split}
& 
\int_{\mathbb R^N} 
|\eta|^{-\varepsilon} 
\sup _{ |y| \leq |\eta|} J_2(y) \frac{d\eta}{|\eta|^N}
\\
& 
\preceq
\int_{\mathbb R^ N} 
\Big( |\eta|^{-\frac{\varepsilon}{p-1}} 
\sup _{ |y| \leq |\eta|} 
\big( \| \triangle_y f \|_{L^\infty} + \| \triangle_y g \|_{L^\infty} 
\big) 
\Big)^{p-1}
\frac{d\eta}{|\eta|^n}
\| f-g \|_{L^q}
\\
& 
\preceq \big( \| f \|_{\dot B^{\frac{\varepsilon}{p-1}}_{\infty,p-1}} ^{p-1} 
            +\| g \|_{\dot B^{\frac{\varepsilon}{p-1}}_{\infty,p-1}} ^{p-1}
       \big)
\| f-g \|_{L^q}
\\
& 
\preceq \Big\{
            \Big( 
             \| f \|_{\dot B^0_{\infty,\infty}} ^{1-\frac{\varepsilon}{p-1}}
             \| f \|_{\dot B^1_{\infty,\infty}} ^{\frac{\varepsilon}{p-1}}
            \Big) ^{p-1}
            +
            \Big( 
             \| g \|_{\dot B^0_{\infty,\infty}} ^{1-\frac{\varepsilon}{p-1}}
             \| g \|_{\dot B^1_{\infty,\infty}} ^{\frac{\varepsilon}{p-1}}
            \Big) ^{p-1}
       \Big\}
\| f-g \|_{L^q}
\\
& 
\preceq \Big(
             \| f \|_{\dot B^0_{\infty, 1}} ^{p-1-\varepsilon}
             \| f \|_{\dot B^1_{\infty, 1}} ^{\varepsilon}
            +
             \| g \|_{\dot B^0_{\infty, 1}} ^{p-1-\varepsilon}
             \| g \|_{\dot B^1_{\infty, 1}} ^{\varepsilon}
       \Big)
\| f-g \|_{\dot B^0_{q,1}} . 
\end{split}
\end{equation}
In the case $ p \geq 2 $, 
by \eqref{eq:2.5} and the H\"older inequality again we obtain 
\begin{equation}\label{eq:2.15}
\begin{split}
\| \triangle_y (|f|^p - |g|^p ) \|_{L^q} 
& 
\preceq 
 \big( \| f \|_{L^\infty} ^{p-1} + \| g \|_{L^\infty} ^{p-1} \big) 
       \| \triangle _y (f-g) \|_{L^q} 
\\
& \quad 
 +  \big( \| f \|_{L^\infty} ^{p-2} + \| g \|_{L^\infty} ^{p-2} \big) 
    \big( \| \triangle _y f \|_{L^\infty} + \| \triangle _y g  \|_{L^\infty} \big) 
   \| f-g \|_{L^q} 
\\
& 
=: J_1(y) + J_3(y) . 
\end{split}
\end{equation}
For the estimate of $J_3$, 
it follows from 
\eqref{eq:2.1} and the embeddings 
$\dot B^{s}_{q,1}  
 \hookrightarrow \dot B^{s}_{q,\infty} $ ($s = 0,1$)
and $\dot B^0_{q,1} \hookrightarrow L^q$ 
 that 
\begin{equation}
\label{eq:2.16}
\begin{split}
& 
\int_{\mathbb R^N}
|\eta|^{-\varepsilon} 
\sup _{ |y| \leq |\eta|} J_3(y) \frac{d\eta}{|\eta|^N}
\\
& 
\preceq \big( \| f \|_{L^\infty} ^{p-2} + \| g \|_{L^\infty} ^{p-2} 
       \big)   
\int_{\mathbb R^N} 
\Big( 
|\eta|^{-\varepsilon} 
\sup _{ |y| \leq |\eta|} 
 \big( \| \triangle _y f \|_{L^\infty} + \| \triangle_y g \|_{L^\infty} 
       \big) 
\Big) 
\frac{d\eta}{|\eta|^N}
       \| f-g \|_{L^q}
\\
& 
\preceq \big( \| f \|_{L^\infty} ^{p-2} + \| g \|_{L^\infty} ^{p-2} 
       \big)   
       \big( \| f \|_{\dot B^\varepsilon_{\infty,1}} + \| g \|_{\dot B^\varepsilon _{\infty,1}} 
       \big) 
       \| f-g \|_{L^q}
\\
& 
\preceq \big( \| f \|_{\dot B^0_{\infty, 1}} ^{p-2} + \| g \|_{\dot B^0_{\infty, 1}} ^{p-2} 
       \big)   
       \big( \| f \|_{\dot B^0_{\infty,1}} ^{1-\varepsilon} 
             \| f \|_{\dot B^1_{\infty,1}} ^\varepsilon
           + \| g \|_{\dot B^0_{\infty,1}} ^{1-\varepsilon} 
             \| g \|_{\dot B^1_{\infty,1}} ^\varepsilon
       \big) 
       \| f-g \|_{\dot B^0_{q, 1}}. 
\end{split}
\end{equation}
Therefore, \eqref{eq:2.4} is obtained by \eqref{eq:2.12}, 
\eqref{eq:2.13}, \eqref{eq:2.14}, \eqref{eq:2.15} and \eqref{eq:2.16}.
$\Box$\vspace{7pt}

The end of this section we recall some results on $e^{t{\mathcal L}}f$.
\begin{Lemma}{\rm\cite{IKK2,Iw}}
\label{Lemma:2.4}
Let $s \in \mathbb R $ and $q \in [1,\infty]$. 

\noindent 
{\rm (i)} 
For $j = 0,1$, $\alpha \geq 0$, $r \in [1,q]$ and $\sigma \in [1,\infty]$, 
it holds that 
\begin{gather}
\label{eq:2.17}
\| \nabla ^j  P_{t+1} \|_{L^q} 
\preceq (1+t)^{- N(1-\frac{1}{q})-j } , 
\\
\label{eq:2.18}
\| \nabla ^j  e^{t\mathcal L} f \|_{L^q} 
\preceq t^{- N(\frac{1}{r}-\frac{1}{q})-j } 
     \| f \|_{L^r} , 
\\ 
\label{eq:2.19}
\| e^{t\mathcal L} f \|_{\dot B^{s+\alpha}_{q,\sigma}} 
\preceq t^{- N(\frac{1}{r}-\frac{1}{q}) - \alpha} 
     \| f \|_{\dot B^s_{r,\sigma}} . 
\end{gather}

\noindent 
{\rm (ii)} 
It holds that 
\begin{equation}\label{eq:2.20}
\| e^{t\mathcal L} f \|_{L^1 _t (0,\infty ; \dot B^{s+1}_{q ,1})} 
\preceq \| f \|_{\dot B^{s}_{q ,1}} . 
\end{equation}

\noindent 
{\rm (iii)} 
It holds that 
\begin{equation}\label{eq:2.21}
\Big\| 
 \int_0^{t} e^{ (t-\tau) \mathcal L} f (\tau) \, d\tau 
\Big\|_{L^1_t (0,\infty ; \dot B^{s+1}_{q ,1})} 
\preceq \| f \|_{L^1 (0,\infty ; \dot B^{s}_{q ,1})}  . 
\end{equation}
\end{Lemma}

\begin{Lemma}
\label{Lemma:2.5}$(${\rm[Proposition~3.1]\cite{IKK}}$)$.
For any $f\in L^1$, if
$$
\int_{{\mathbb R}^N}f(x)\,dx=0,
$$
then
$$
\lim_{t\to\infty}\|e^{t{\cal L}}f\|_{L^1}=0.
$$
\end{Lemma}
%

\section{Existence of global-in-time solutions and Decay estimates}

In this section we prove the assertion {\rm (i)} of Theorem \ref{Theorem:1.1}. 
We apply the contraction mapping principle in 
a suitable complete metric space. Let $\Psi (u)$ be defined by 
\begin{equation}
\label{eq:3.1}
\Psi (u) (t) 
:= e^{t\mathcal L} u_0 + 
  \int_0^t e^{(t-\tau)\mathcal L} |\nabla u(\tau)|^p \,d\tau , 
\end{equation}
and we define the following norms  
\begin{equation}\notag 
\begin{split}
\| u \|_{\dot X^s_q} 
& 
:= \sup _{t > 0} \| u(t) \|_{\dot B^s_{q,1}} 
  +   \int_0^\infty \| u(t) \|_{\dot B^{s+1}_{q,1}} dt ,
\\
\| u \|_{\dot Y^s_q} 
& :=  \| u \|_{\dot X^s_q}
  + \sup _{t > 0} t^{N(\frac{1}{r}-\frac{1}{q})+1-s} \| u(t) \|_{\dot B^1_{q,1}}
  + \int_0^\infty t^{N(\frac{1}{r}-\frac{1}{q})+1-s} \| u(t) \|_{\dot B^{2}_{q,1}} dt .
\end{split}
\end{equation}
$\dot X^s_q$ with the norm $\| u \|_{\dot X^s_q}$ 
is defined by the space of all functions $u$ such that  
$$
u \in L^\infty(0,\infty ; \dot B^s_{q,1} ) \cap 
     L^1 (0,\infty ; \dot B^{s+1}_{q,1}) 
   \quad  \text{and} \quad 
 \| u \|_{\dot X^s_q} < \infty, 
$$
and 
$\dot Y^s_q$ with the norm $\| \cdot \|_{\dot Y^s_q}$ 
is also defined by the space of all functions $u$ such that 
$$
u \in \dot X^s_q
\quad \text{and} \quad 
\| u \|_{\dot Y^s_q} < \infty. 
$$
Here, let $\varepsilon$ and $\lambda$ be fixed constants satisfying  
\begin{equation}
\label{eq:3.2}
0< \varepsilon < \min\{1,p-1 \} 
\quad \text{and} \quad 0 < \lambda <1, 
\end{equation}
and we introduce 
the following metric space $\mathfrak X$
\begin{gather}\notag 
\mathfrak X 
:= \big\{ u \in \dot X^\varepsilon_r \cap \dot X^\varepsilon_\infty 
               \cap \dot X^1_r \cap \dot X^ 1_\infty  
      \, \big| \, 
      \| u \|_{\dot X^1_q } \leq 2 C_0 \| u_0 \|_{\dot B^1_{q,1}}  
       \text{ for any } q \in [r,\infty] , 
\\ \notag 
\hspace*{44mm}      \| u \|_{\dot Y^\lambda_r \cap \dot Y^\lambda _{\infty}} 
       \leq 2 C_0 \| u_0 \|_{\dot B^\lambda_{r,1} \cap \dot B^\lambda _{\infty,1} } 
    \big\}, 
\end{gather}
with the metric 
\begin{gather} \notag 
 d(u,v) := \| u-v \|_{\dot X^\varepsilon_r \cap \dot X^\varepsilon_\infty}, 
\end{gather}
where $C_0$ will be taken later. 
We first show that $\mathfrak X$ is a complete metric space. 

\begin{Lemma}
$\mathfrak X$ is a complete metric space. 
\end{Lemma}

\noindent 
{\bf Proof. } 
It is easy to see that $\mathfrak X$ is a metric space, 
then we show the completeness only. 
Let $\{ u_n \}$ be a Cauchy sequence in $\mathfrak X$. 
Since $\dot X^\varepsilon_r \cap \dot X^\varepsilon_\infty$ is complete, 
there exists $u \in \dot X^\varepsilon_r \cap \dot X^\varepsilon_\infty$ 
such that $u_n $ converges to $u$ in $\dot X^\varepsilon_r \cap \dot X^\varepsilon_\infty$ 
as $n \to \infty$. Then we also have 
\begin{gather}\notag 
\lim _{n\to \infty} \| \phi_j * ( u_n(t) - u(t) ) \|_{L^r \cap L^\infty} 
 = 0 
\quad \text{for almost every } t \text{ and } j \in \mathbb Z, 
\\ \notag 
\lim_{n\to \infty} 
\int_0^L \| \phi_j * ( u_n(t) - u(t) ) \|_{L^r \cap L^\infty}\,dt 
= 0 
\quad \text{for any }  L > 0 \text{ and } j \in \mathbb Z.
\end{gather}
There imply that,
for any $L > 0$, 
\begin{eqnarray*} 
&&
\lim _{n \to \infty}
\sum _{|j| \leq L} 2^j \| \phi_j * u_n(t) \|_{L^r \cap L^\infty}
= 
\sum _{|j| \leq L} 2^j \| \phi_j * u(t) \|_{L^r \cap L^\infty} , 
\\
&& 
\lim_{n\to \infty} 
t^{N(\frac{1}{r}-\frac{1}{q})+1-\lambda} 
\sum _{|j| \leq L} 2^j \| \phi_j * u_n(t) \|_{L^ q}
= 
t^{N(\frac{1}{r}-\frac{1}{q})+1-\lambda} 
\sum _{|j| \leq L} 2^j \| \phi_j * u_n(t) \|_{L^ q} 
\end{eqnarray*}
for almost every $t$, and 
\begin{eqnarray*}
&&
\lim_{n\to \infty} 
\int_0^L \sum_{|j| \leq L} 2^{2j}\| \phi_j * u_n(t) \|_{L^r \cap L^\infty}\,dt 
= 
\int_0^L \sum_{|j| \leq L} 2^{2j}\| \phi_j * u(t) \|_{L^r \cap L^\infty}\,dt , 
\\ 
&&
\begin{aligned}
\lim_{n\to \infty} 
& 
\int_0^L t^{N(\frac{1}{r}-\frac{1}{q})+1-\lambda}  
  \sum_{|j| \leq L} 2^{2j}\| \phi_j * u_n(t) \|_{L^r \cap L^\infty}\,dt 
\\
&\qquad\qquad\qquad\qquad
= 
\int_0^L 
  t^{N(\frac{1}{r}-\frac{1}{q})+1-\lambda}  
   \sum_{|j| \leq L} 2^{2j}\| \phi_j * u(t) \|_{L^r \cap L^\infty}\,dt . 
\end{aligned}
\end{eqnarray*}
The terms in right hand side of the above four equalities are bounded uniformly 
with respect to $L > 0$ since $\{ u_n \} \subset \mathfrak X$, 
and they are monotone increasing, so that, they converges as $L \to \infty$. 
Then we deduce that $u$ satisfies 
$$
\| u \|_{\dot X^1_q } \leq 2 C_0 \| u_0 \|_{\dot B^1_{q,1}}  
       \text{ for any } q \in [r,\infty] 
\quad \text{and} \quad     
\| u \|_{\dot Y^\lambda_r \cap \dot Y^\lambda _{\infty}} 
       \leq 2 C_0 \| u_0 \|_{\dot B^\lambda_{r,1} \cap \dot B^\lambda _{\infty,1} } , 
$$
hence, $u \in \mathfrak X$. 
Therefore the completeness of $\mathfrak X$ follows. 
$\Box$\vspace{7pt}

In order to estimate the terms in \eqref{eq:3.1}, 
we prepare the following proposition.

\begin{Proposition}\label{Proposition:3.1}
Let $p$, $q$, $r$, $\varepsilon$ and $\lambda$ satisfy 
$p > 1$,  $1 \leq r \leq q \leq \infty$ and \eqref{eq:3.2}.  
Then there holds that 
\begin{eqnarray}
&&
\label{eq:3.3}
\| e^{t\mathcal L} u_0 \|_{\dot X^1_q} 
\preceq \| u_0 \|_{\dot B^1_{q,1}}, 
\qquad\qquad
\| e^{t\mathcal L} u_0 \|_{\dot Y^\lambda_q} 
\preceq \| u_0 \|_{\dot B^\lambda_{q,1}}+\| u_0 \|_{\dot B^\lambda_{r,1}}, 
\\ 
&&
\label{eq:3.4}
\Big\| 
  \int_0^t e^{(t-\tau)\mathcal L} |\nabla u(\tau)|^p \,d\tau
\Big\|_{L^q} 
\preceq t \| u \|_{\dot X^1_{\infty}}^{p-1} \| u \|_{\dot X^1_{q}}, 
\\ 
&&
\label{eq:3.5}
\Big\| 
  \int_0^t e^{(t-\tau)\mathcal L} |\nabla u(\tau)|^p \,d\tau 
\Big\|_{\dot X^1_q} 
\preceq \| u \|_{\dot X^1_{\infty}} ^{p-1}\| u \|_{\dot X^1_{q}} , 
\\ 
&&
\label{eq:3.6}
\Big\| 
  \int_0^t e^{(t-\tau)\mathcal L} |\nabla u(\tau)|^p \,d\tau 
\Big\|_{\dot Y^\lambda_q} 
\preceq \| u \|_{\dot X^1_{\infty}} ^{p-1}
   \big( \| u \|_{\dot Y^\lambda_r} + \| u \|_{\dot Y^\lambda_q} \big) , 
\\ 
&&
\label{eq:3.7}
\Big\| 
  \int_0^t e^{(t-\tau)\mathcal L} \big( |\nabla u(\tau)|^p - |\nabla v(\tau)|^p \big) \,d\tau
\Big\|_{\dot X^\varepsilon_q} 
\preceq
   \big( \| u \|_{\dot X^1_\infty} ^{p-1} + \| v \|_{\dot X^1_{\infty}} ^{p-1}
   \big)
   \| u-v \|_{\dot X^\varepsilon_q} .  
\end{eqnarray}
\end{Proposition}
\begin{Remark}
We should note that the nonlinear term $|\nabla u(\tau)|^p$ 
is in $\mathcal S' (\mathbb R^N)$ if $u \in \dot X^1_\infty$. 
Although $\dot B^1_{\infty,1} (\mathbb R^N)$ is considered as 
a subspace of $\mathcal Z'(\mathbb R^N)$, 
$\nabla u (\tau)$ is determined 
independently of the choice of 
representative elements in $\dot B^1_{\infty,1}$ 
by $\nabla u(\tau) \in \dot B^0_{\infty,1} (\mathbb R^N)$ and Remark $\ref{Remark:2.1}$, 
hence, $\nabla u(\tau) \in \mathcal S '(\mathbb R^N)$. 
We also see  $|\nabla u(\tau)|^p \in L^\infty (\mathbb R^N) 
\subset \mathcal S'(\mathbb R^N)$ 
by $\dot B^0_{\infty,1} \subset L^\infty$. 
In addition the estimate \eqref{eq:3.4} implies that 
the term in the left hand side is in $L^q$ if $u \in \mathfrak X$. 
\end{Remark}
\noindent
{\bf Proof.} 
The linear estimate \eqref{eq:3.3} is verified by the use of 
\eqref{eq:2.19} and \eqref{eq:2.20}. 
In fact, the estimates of $\| e^{t\mathcal L} u_0 \|_{\dot X^{s}_q}$ $(s = 1,\lambda)$ 
is obtained by the boundedness of $e^{t\mathcal L}$, \eqref{eq:2.19} and \eqref{eq:2.20}, 
and for the second and third terms in the definition of $\| \cdot \|_{\dot Y^\lambda _q}$ 
we also apply \eqref{eq:2.19} and \eqref{eq:2.20} to get 
\begin{eqnarray*}
&&
t^{N(\frac{1}{r}-\frac{1}{q})+1-\lambda} \| e^{t\mathcal L} u_0 \|_{\dot B^1_{q,1}} 
\preceq \| u_0 \|_{\dot B^\lambda_{r,1}}, 
\\
&&
\begin{aligned}
\int_0^\infty t^{N(\frac{1}{r}-\frac{1}{q})+1-\lambda} \| e^{t\mathcal L} u_0 \|_{\dot B^{2}_{q,1}} \, dt 
& 
\preceq \int_0^\infty \| e^{\frac{t}{2} \mathcal L } u_0 \|_{\dot B^{1+\lambda}_{r,1}} \, dt
\preceq \| u_0 \|_{\dot B^\lambda_{r,1}} . 
\end{aligned}
\end{eqnarray*}
Then \eqref{eq:3.3} is obtained. 

The estimate \eqref{eq:3.4} is obtained by applying the boundedness of 
$e^{(t-\tau)\mathcal L} $ from $L^q$ to itself, 
the H\"older inequality and the embedding 
$\dot B^0_{q,1} \hookrightarrow L^q$. 
Thus we omit the detail. 

To prove the nonlinear estimates \eqref{eq:3.5}, \eqref{eq:3.6} and \eqref{eq:3.7}, 
we prepare the following nonlinear estimates that, for $s = 1,\lambda $,
\begin{eqnarray}
&&
\label{eq:3.8}
\| |\nabla u|^p \|_{\dot B^{s}_{ q,1}}
\preceq 
\| \nabla u \|_{\dot B^0_{\infty,1}}^{p-1} \| \nabla u \|_{\dot B^{s}_{q,1}} 
\preceq 
\| u \|_{\dot B^1_{\infty,1}}^{p-1} \| u \|_{\dot B^{s +1}_{q,1}},
\\ 
&&
\label{eq:3.9}
\begin{aligned}
& \big\| |\nabla u|^p - |\nabla v|^p 
 \big\|_{\dot B^\varepsilon_{ q,1}}
\\
& 
\preceq 
       \big( \| \nabla u\|_{\dot B^0_{\infty,1}} ^{p-1} + \| \nabla v \|_{\dot B^0_{\infty,1}}^{p-1} 
       \big) 
       \| \nabla u -\nabla v \|_{\dot B^\varepsilon_{q,1}}
\\
& \quad 
  + \Big(
             \| \nabla u \|_{\dot B^0_{\infty, 1}} ^{p-1-\varepsilon}
             \| \nabla u \|_{\dot B^1_{\infty, 1}} ^{\varepsilon}
            +
             \| \nabla v \|_{\dot B^0_{\infty, 1}} ^{p-1-\varepsilon}
             \| \nabla v \|_{\dot B^1_{\infty, 1}} ^{\varepsilon}
       \Big)
\| \nabla u - \nabla v \|_{\dot B^0_{q,1}}
\\
& 
\preceq 
       \big( \| u \|_{\dot B^1_{\infty,1}} ^{p-1} + \| v \|_{\dot B^1_{\infty,1}}^{p-1} 
       \big) 
       \| u - v \|_{\dot B^{1+\varepsilon}_{q,1}}
\\
& \quad 
+\| u - v \|_{\dot B^{\varepsilon}_{q,1}}  ^{\varepsilon}
\| u - v \|_{\dot B^{1+\varepsilon}_{q,1}} ^{1-\varepsilon}
\times
\\
&
\qquad
\times
\left\{
\begin{array}{l}
\displaystyle{
\Big(
             \| u \|_{\dot B^1_{\infty, 1}} ^{p-1-\varepsilon}
             \| u \|_{\dot B^2_{\infty, 1}} ^{\varepsilon}
            +
             \| v \|_{\dot B^1_{\infty, 1}} ^{p-1-\varepsilon}
             \| v \|_{\dot B^2_{\infty, 1}} ^{\varepsilon}
       \Big)},
\hspace{2.6cm} 
\text{if}\quad 1 < p < 2, 
\vspace{5pt}\\
\displaystyle{
\Big( \| u \|_{\dot B^1_{\infty, 1}} ^{p-2} + \| v \|_{\dot B^1_{\infty, 1}} ^{p-2} 
       \Big)   
       \Big( \| u \|_{\dot B^1_{\infty,1}} ^{1-\varepsilon} 
             \| u \|_{\dot B^2_{\infty,1}} ^\varepsilon
           + \| v \|_{\dot B^1_{\infty,1}} ^{1-\varepsilon} 
             \| v \|_{\dot B^2_{\infty,1}} ^\varepsilon
       \Big) },
       \,\,
\text{if}\quad p \geq 2 , 
\end{array}
\right. 
\end{aligned}
\end{eqnarray}
which are obtained by \eqref{eq:2.3}, \eqref{eq:2.4} and the interpolation inequality 
in the Besov spaces, that is, 
$$
\|f\|_{\dot B^1_{q,1}}=
\sum _{ j \in \mathbb Z}2^j  \| \phi_j * f \|_{L^q}
= \sum _{ j \in \mathbb Z} ( 2^{\varepsilon j} ) ^\varepsilon (2^{(1+\varepsilon) j} )^{1-\varepsilon} 
  \| \phi_j * f \|_{L^q}^{\varepsilon} \| \phi_j * f \|_{L^q} ^{1-\varepsilon}
\leq 
\|f\|_{\dot B^\varepsilon_{q,1}} ^\varepsilon 
\| f \|_{\dot B^{1+\varepsilon}_{q,1}} ^{1-\varepsilon}. 
$$
  
On the estimate of \eqref{eq:3.5},
by the boundedness of $e^{(t-\tau) \mathcal L}$, \eqref{eq:2.21}, 
\eqref{eq:3.8} with $s = 1$ and the H\"older inequality 
we have
\begin{equation}\notag 
\begin{split}
\Big\| 
  \int_0^t e^{(t-\tau)\mathcal L} |\nabla u(\tau)|^p  \, d\tau
\Big\|_{\dot X^1_q} 
& 
\preceq \int_0^\infty \| |\nabla u(\tau)|^p \|_{\dot B^1_{ q,1}} \, d\tau
\\
& 
\preceq \| u \|_{L^\infty (0,\infty;\dot B^1_{\infty,1})} ^{p-1} 
       \| u \|_{L^1 (0,\infty ; \dot B^2_{ q,1})}
\\
& 
\preceq \| u \|_{\dot X^1_{\infty}} ^{p-1}\| u \|_{\dot X^1_{q}} . 
\end{split}
\end{equation}
Then \eqref{eq:3.5} is obtained. 

On the estimate of \eqref{eq:3.6}, the first norm $\| \cdot \|_{\dot X^\lambda _q}$ 
in the definition of $\| \cdot \|_{\dot Y^\lambda _q}$
can be treated in the same way as the proof of \eqref{eq:3.5} 
with \eqref{eq:3.8} ($s = \lambda$), 
thus we omit the estimate on $\| \cdot  \|_{\dot X^\lambda _q}$ 
to consider the second and third terms only. 
We put
\begin{equation}
\label{eq:3.10}
K_1(x,t):=\int_0^{t/2} e^{(t-\tau)\mathcal L} |\nabla u(\tau)|^p  \, d\tau,
\qquad
K_2(x,t):=\int_{t/2}^t e^{(t-\tau)\mathcal L} |\nabla u(\tau)|^p  \, d\tau,
\end{equation}
for all $x\in{\mathbb R}^N$ and $t>0$.
On the second term, 
by \eqref{eq:2.19} and \eqref{eq:3.8} we have
\begin{equation}
\label{eq:3.11}
\begin{split}
&
t^{N(\frac{1}{r}-\frac{1}{q}) + 1 -\lambda}
\| K_1(t)\|_{\dot B^1_{q,1}} 
\\
& 
\preceq 
t^{N(\frac{1}{r}-\frac{1}{q}) + 1 -\lambda}
\int_0^{t/2} (t-\tau)^{-N(\frac{1}{r}-\frac{1}{q}) -1 + \lambda}
\big\| |\nabla u (\tau)|^p  \big\|_{\dot B^\lambda_{r,1}} \, d\tau
\\
& 
\preceq \int_0^{t/2} \| u(\tau) \|_{\dot B^1_{\infty,1}} ^{p-1} \| u(\tau) \|_{\dot B^{\lambda + 1}_{r,1}} 
\, d\tau
\\
& \preceq \| u \|_{\dot X^1_{\infty}} ^{p-1}\| u \|_{\dot Y^\lambda_{r}},\qquad t>0.
\end{split}
\end{equation}
Furthermore,
by the boundedness of $e^{(t-\tau)\mathcal L}$ and \eqref{eq:3.8}
we obtain
\begin{equation}
\label{eq:3.12}
\begin{split}
&
t^{N(\frac{1}{r}-\frac{1}{q}) + 1 -\lambda}
\| K_2(t)\|_{\dot B^1_{q,1}} 
\\
& 
\preceq t^{N(\frac{1}{r}-\frac{1}{q}) + 1 -\lambda}
\int_{t/2}^t 
\big\| |\nabla u (\tau)|^p  \big\|_{\dot B^1_{q,1}} \, d\tau
\\
& 
\preceq \| u \|_{L^\infty (0,\infty ; \dot B^1_{\infty,1})} 
       \int_{t/2}^t \tau^{N(\frac{1}{r}-\frac{1}{q}) + 1 -\lambda}
       \| u(\tau) \|_{\dot B^2_{q,1}} 
       \,d\tau
\\
& 
\preceq \| u\|_{\dot X^1_\infty} ^{p-1} \| u \|_{\dot Y^\lambda_q},\qquad t>0 . 
\end{split}
\end{equation}
On the third norm, by \eqref{eq:2.19}, \eqref{eq:2.20} and\eqref{eq:3.8} we have 
\begin{equation}
\label{eq:3.13}
\begin{split}
&
\int_0^\infty t^{N(\frac{1}{r}-\frac{1}{q}) + 1 -\lambda} 
\|K_1(t)\|_{\dot B^2_{q,1}} 
dt
\\
& 
\preceq 
\int_0^\infty t^{N(\frac{1}{r}-\frac{1}{q}) + 1 -\lambda}
\int_0^{t/2} (t-\tau)^{-N(\frac{1}{r}-\frac{1}{q}) -1 + \lambda}
\big\| e^{\frac{t-\tau}{2}\mathcal L}|\nabla u (\tau)|^p  \big\|_{\dot B^{\lambda + 1}_{r,1}} \, d\tau \,dt
\\
& 
\preceq 
\int_0^\infty 
\int_{2\tau}^{\infty} 
\big\| e^{\frac{t-\tau}{2}\mathcal L}|\nabla u (\tau)|^p  \big\|_{\dot B^{\lambda + 1}_{r,1}} \, dt\, d\tau
\preceq \int_0^\infty  \| | \nabla u (\tau) |^p \|_{\dot B^{\lambda}_{r,1}} 
\, d\tau
\\
& \preceq \| u \|_{\dot X^1_{\infty}} ^{p-1}\| u \|_{\dot Y^\lambda_{r}}.
\end{split}
\end{equation}
Furthermore, by \eqref{eq:2.20} and \eqref{eq:3.8} we obtain
\begin{equation}
\label{eq:3.14}
\begin{split}
& 
\int_0^\infty t^{N(\frac{1}{r}-\frac{1}{q}) + 1 -\lambda}
\|K_2(t) \|_{\dot B^2_{q,1}} dt 
\\
& 
\preceq 
\int_0^\infty \tau^{N(\frac{1}{r}-\frac{1}{q}) + 1 -\lambda} \int_{\tau}^{2\tau} 
\big\| e^{(t-\tau)\mathcal L}|\nabla u (\tau)|^p  \big\|_{\dot B^2_{q,1}} \, dt \, d\tau
\\
& 
\preceq 
\int_0^\infty \tau^{N(\frac{1}{r}-\frac{1}{q}) + 1 -\lambda}
\big\| |\nabla u (\tau)|^p  \big\|_{\dot B^1_{q,1}} \, d\tau
\\
& 
\preceq 
\| u \|_{L^\infty (0,\infty ; \dot B^1_{\infty,1})} ^{p-1}
\int_0^\infty \tau^{N(\frac{1}{r}-\frac{1}{q}) + 1 -\lambda}
\| u(\tau) \|_{\dot B^2_{q,1}} \, d\tau
\\
& 
\preceq 
\| u \|_{\dot X^1_{\infty}} ^{p-1} \| u \|_{\dot Y^\lambda _q} . 
\end{split}
\end{equation}
Then, \eqref{eq:3.6} is obtained by \eqref{eq:3.10}, \eqref{eq:3.11}, \eqref{eq:3.12}, \eqref{eq:3.13} and \eqref{eq:3.14}.

For the proof of \eqref{eq:3.7}, 
it follows from 
the boundedness of $e^{(t-\tau)\mathcal L}$, \eqref{eq:2.21}, \eqref{eq:3.9} 
and the H\"older inequality that, 
if $1<p<2$, then
\begin{equation}\notag 
\begin{split}
& 
\Big\| 
  \int_0^t e^{(t-\tau)\mathcal L} \Big( |\nabla u(\tau)|^p - |\nabla v(\tau)|^p \Big) \, d\tau
\Big\|_{\dot X^\varepsilon_q} 
\\
& 
\preceq \int_0^\infty \big\| |\nabla u(\tau)|^p - |\nabla v(\tau)|^p \big\|_{\dot B^\varepsilon_{ q,1}} \, d\tau
\\
& 
\preceq \Big( \| u \|_{L^\infty (0,\infty ; \dot B^1_{\infty,1})} ^{p-1}
            +\| v \|_{L^\infty (0,\infty ; \dot B^1_{\infty,1})} ^{p-1}
       \Big) 
       \| u-v \|_{L^1(0,\infty ; \dot B^{1+\varepsilon}_{q,1})} 
\\
& 
\quad + \| u-v \|_{L^\infty (0,\infty ; \dot B^{\varepsilon}_{q,1})} ^{\varepsilon}
        \| u-v \|_{L^1 (0,\infty ; \dot B^{1+\varepsilon}_{q,1})} ^ {1-\varepsilon}
        \times
\\
& \qquad
\times
 \Big( \| u \|_{L^\infty (0,\infty ; \dot B^1_{\infty,1})} ^{p-1-\varepsilon}
                \| u \|_{L^1(0,\infty ; \dot B^2_{\infty,1})} ^{\varepsilon}
               +\| v \|_{L^\infty (0,\infty ; \dot B^1_{\infty,1})} ^{p-1-\varepsilon}
                \| v \|_{L^1(0,\infty ; \dot B^2_{\infty,1})} ^{\varepsilon}
          \Big)
\\
& 
\preceq \Big( \| u \|_{\dot X^1_{\infty}} ^{p-1}
            +\| v \|_{\dot X^1_{\infty}} ^{p-1}
       \Big) 
       \| u-v \|_{\dot X^{\varepsilon}_{ q}},
\end{split}
\end{equation}
and, if $p\ge2$, then
\begin{equation}\notag 
\begin{split}
& 
\Big\| 
  \int_0^t e^{(t-\tau)\mathcal L} \Big( |\nabla u(\tau)|^p - |\nabla v(\tau)|^p \Big) \, d\tau
\Big\|_{\dot X^\varepsilon_q} 
\\
& 
\preceq \Big( \| u \|_{L^\infty (0,\infty ; \dot B^1_{\infty,1})} ^{p-1}
            +\| v \|_{L^\infty (0,\infty ; \dot B^1_{\infty,1})} ^{p-1}
       \Big) 
       \| u-v \|_{L^1(0,\infty ; \dot B^{1+\varepsilon}_{q,1})} 
\\
& 
\quad +  \| u-v \|_{L^\infty (0,\infty ; \dot B^{\varepsilon}_{q,1})} ^{\varepsilon}
        \| u-v \|_{L^1 (0,\infty ; \dot B^{1+\varepsilon}_{q,1})} ^ {1-\varepsilon}
\Big( \| u \|_{L^\infty (0,\infty ; \dot B^1_{\infty,1})} ^{p-2}+\| v \|_{L^\infty (0,\infty ; \dot B^1_{\infty,1})} ^{p-2}
	\Big)
	\times
\\
& \qquad
              \times  \Big(\| u \|_{L^\infty (0,\infty ; \dot B^1_{\infty,1})} ^{1-\varepsilon}\| u \|_{L^1(0,\infty ; \dot B^2_{\infty,1})} ^{\varepsilon}
                +\| v \|_{L^\infty (0,\infty ; \dot B^1_{\infty,1})} ^{1-\varepsilon}\| v \|_{L^1(0,\infty ; \dot B^2_{\infty,1})} ^{\varepsilon}
                \Big)
\\
& 
\preceq \Big( \| u \|_{\dot X^1_{\infty}} ^{p-1}
            +\| v \|_{\dot X^1_{\infty}} ^{p-1}
       \Big) 
       \| u-v \|_{\dot X^{\varepsilon}_{ q}},
\end{split}
\end{equation}
Therefore, \eqref{eq:3.7} is obtained and the proof of all estimates is completed. 
$\Box$\vspace{7pt}

In what follows, 
we prove that the solution exists globally in time by applying 
the contraction mapping principle in $\mathfrak X$ for initia 
data $u_0$ in $B^1_{r,1} \cap B^1_{\infty,1} $ and small in $\dot B^1_{\infty,1}$, 
and that the solutions satisfy the decay estimates \eqref{eq:1.9} and \eqref{eq:1.10}. 

\vspace{7pt}

\noindent 
{\bf Proof of existence of global-in-time solutions in $\mathfrak X 
\cap C([0,\infty), B^1_{r,1} \cap B^1_{\infty,1})$. } 
Let the constant $C_0$ in the definition of $\mathfrak X$ 
be a constant which satisfy the all estimates in Proposition \ref{Proposition:3.1}, 
and we assume the initial data satisfies 
\begin{equation}\label{eq:3.15}
u_0 \in B^1_{r,1} \cap B^1_{\infty,1} 
\quad \text{and} \quad 
\| u_0 \|_{\dot B^1_{\infty,1 }} 
\leq (2^{p+1} C_0 ^p) ^{-\frac{1}{p-1}} . 
\end{equation}
For any $u , v \in \mathfrak X$, it follows from Proposition \ref{Proposition:3.1} that 
\begin{eqnarray}\notag
&& 
\begin{aligned}
\| \Psi(u) \|_{\dot X^s_q} 
& \leq C_0 \| u_0 \|_{\dot B^s_{q,1}} 
  + C_0 \| u \|_{\dot X^1_\infty} ^{p-1} \| u \|_{\dot X^s_q}
\\
& \leq C_0 \| u_0 \|_{\dot B^s_{q,1}} 
  + C_0 (2C_0 \| u_0 \|_{\dot B^1_{\infty,1}}) ^{p-1} \cdot 2C_0 \| u_0 \|_{\dot B^s_{q,1}}
\\
& \leq 2 C_0 \| u_0 \|_{\dot B^s_{q,1}} , 
\end{aligned}
\\ \notag 
&&
\begin{aligned}
\| \Psi(u) \|_{\dot  Y^\lambda_r \cap \dot Y^\lambda _\infty} 
& \leq C_0 \| u_0 \|_{\dot B^\lambda_{r,1} \cap \dot B^\lambda_{\infty,1}} 
  + C_0 \| u \|_{\dot X^1_\infty} ^{p-1} \| u \|_{\dot Y^\lambda_{r}}
\\
& \leq C_0 \| u_0 \|_{\dot B^\lambda_{r,1} \cap \dot B^\lambda_{\infty,1}} 
  + C_0 (2C_0 \| u_0 \|_{\dot B^1_{\infty,1}}) ^{p-1} \cdot 2C_0 
    \| u_0 \|_{\dot B^\lambda_{r,1} \cap \dot B^\lambda_{\infty,1}} 
\\
& \leq 2 C_0 \| u_0 \|_{\dot B^\lambda_{r,1} \cap \dot B^\lambda_{\infty,1}} , 
\end{aligned}
\\
&&
 \label{eq:3.16}
\begin{aligned}
\| \Psi(u) - \Psi(v)\|_{\dot X^\varepsilon_r \cap \dot X^\varepsilon_\infty}  
& \leq C_0 (\| u \|_{\dot X^1_{\infty} } ^{p-1} + \| v \|_{\dot X^1_\infty} ^{p-1}) 
    \| u-v \|_{\dot X^\varepsilon_r \cap \dot X^\varepsilon_\infty}  
\\
& \leq C_0 \cdot 2 (2C_0 \| u_0 \|_{\dot B^1_{\infty,1}}) ^{p-1}  \cdot 
    \| u-v \|_{\dot X^\varepsilon_r \cap \dot X^\varepsilon_\infty}  
\\& 
\leq \frac{1}{2} \| u-v \|_{\dot X^\varepsilon_r \cap \dot X^\varepsilon_\infty} . 
\end{aligned}
\end{eqnarray}
for any $s = \varepsilon, 1$.
$\Psi$ is a contraction map from $\mathfrak X$ to itself 
and the global solution for small initial data is obtained in $\mathfrak X$. 
Then $u(t) = \Psi (u)(t)$ in $\mathcal Z'(\mathbb R^N)$ for allmost 
every $t$, and 
we have to find a fixed point such that the equality $u(t) = \Psi (u)(t)$ 
holds in $\mathcal S'(\mathbb R^N)$. 
For this purpose, we take a sequence $\{ u_n \}$ such that 
$$
u_1 := e^{t\mathcal L} u_0 , 
\quad 
u_n := \Psi(u_{n-1}), 
  \quad n \geq 2.
$$
The previous contraction argument implies that 
$u_n$ converges to $u$ in $\dot X^\varepsilon_r \cap \dot X^\varepsilon_\infty$. 
Here, we see that $\Psi (u_{n-1}) (t)$ tends to $\Psi (u) (t)$ in $L^\infty$ 
as $n \to \infty$ for each $t$ since we have from \eqref{eq:3.4}
\begin{eqnarray*}
&& 
\| \Psi (u_{n-1}) (t) \|_{L^\infty} 
\preceq \| u_0 \|_{L^\infty} + t \| u_{n-1} \|_{\dot X^1_{\infty}} ^{p} , 
\quad 
\| \Psi (u) (t) \|_{L^\infty} 
\preceq \| u_0 \|_{L^\infty} + t \| u \|_{\dot X^1_{\infty}} ^{p} , 
\\
&&
\begin{aligned}
& 
\Big\| 
\int_0^t e^{(t-\tau)\mathcal L} 
  ( |\nabla u_{n-1}(\tau)|^p - |\nabla u (\tau)|^p ) d\tau 
\Big\|  _{L^\infty}
\\
& 
\preceq 
\int_0^t 
  \Big( \| \nabla u_{n-1} \|_{L^\infty} ^{p-1} + \| \nabla u(\tau) \|_{L^\infty} ^{p-1}
  \Big) 
  \| \nabla u_{n-1} (\tau) - \nabla u(\tau) \|_{L^\infty}
  d\tau
\\
& 
\preceq \Big( \| u_{n-1} \|_{\dot X^1_{\infty}}^{p-1} + \| u \|_{\dot X^1_{\infty}}^{p-1}
        \Big)
\int_0^t 
  \| u_{n-1} (\tau) - u(\tau) \|_{\dot B^1_{\infty,1}} d\tau
\\
& 
\preceq \| u_0 \|_{\dot B^1_{\infty,1}} ^{p-1} 
\int_0^t 
  \| u_{n-1} (\tau) - u(\tau) \|_{\dot B^\varepsilon_{\infty,1}}  ^\varepsilon
  \| u_{n-1} (\tau) - u(\tau) \|_{\dot B^{1+\varepsilon}_{\infty,1}} ^{1-\varepsilon}
  d\tau 
\\
& 
\preceq \| u_0 \|_{\dot B^1_{\infty,1}} ^{p-1} 
t^\varepsilon 
\| u_{n-1} - u \|_{L^\infty (0,\infty ; \dot B^\varepsilon_{\infty,1})} ^\varepsilon
\| u_{n-1} - u \|_{L^1 (0,\infty ; \dot B^{1+\varepsilon}_{\infty,1})} ^{1-\varepsilon}
\\
& 
\preceq \| u_0 \|_{\dot B^1_{\infty,1}} ^{p-1} 
t^\varepsilon 
\| u_{n-1} - u \|_{\dot X^\varepsilon _\infty} 
\to 0 
\quad \text{as } n \to \infty . 
\end{aligned}
\end{eqnarray*}
Therefore, $u_n (t)= \Psi (u_{n-1}) (t)$ is a Cauchy sequence in $L^\infty$, 
so that, there exists $v(t) \in L^\infty$ such that 
$u_n (t)$ converges to $v(t)$ in $L^\infty$ as $n \to \infty$. 
It follows from $L^\infty \subset \mathcal S'(\mathbb R^N ) \subset \mathcal Z '(\mathbb R^N)$ 
and the uniqueness of the limit in $\mathcal Z'(\mathbb R^N)$
that 
$u_n (t)$ also converges to $v(t)$ in $\mathcal Z'(\mathbb R^N)$ as $n \to \infty$ 
and $v(t) = u(t)$ in $\mathcal Z'(\mathbb R^N)$. 
Since $u(t) \in \dot B^1_{\infty,1}$ and 
$\nabla u(t) \in \mathcal S'(\mathbb R^N)$ by Remark~\ref{Remark:2.1}, 
it holds that $\nabla v (t) = \nabla u(t)$ 
and $\Psi (v) (t) = \Psi (u) (t)$ in $\mathcal S'(\mathbb R^N)$ for all $t$. 
Then taking the limit in the topology of $L^\infty$ on the equation 
$u_n(t) = \Psi (u_{n-1}) (t)$ for each $t$, we obtain 
$$
v(t) = \Psi (u)(t) = \Psi (v) (t) 
\quad \text{in } L^\infty.
$$
By $\| e^{t\mathcal L} u_0  \|_{L^q} \leq \| u_0 \|_{L^q}$ and \eqref{eq:3.4}, 
$v$ satisfies $v(t) \in 
L^r \cap \dot B^1_{r,1} \cap L^\infty \cap \dot B^1_{\infty,1} 
= B^1_{r,1} \cap B^1_{\infty,1}$. 
Hence, the fixed point $v$ is a solution in $\mathfrak X \cap 
C([0,\infty) , B^1_{r,1} \cap B^1_{\infty,1})$.

It remains to show the uniqueness. 
Let $u, v \in \mathfrak X$ satisfies 
$u, v\in C([0,\infty), B^1_{r,1} \cap B^1_{\infty,1})$, 
$u = \Psi (u)$ and $v = \Psi (v)$, 
and we show that $u(t) = v(t)$ in $\mathcal S'(\mathbb R^N)$ for all $t$. 
The contraction property \eqref{eq:3.16} implies that 
$u(t) = v(t) $ in $\dot B^\varepsilon_{\infty,1} \subset \mathcal Z'(\mathbb R^N)$. 
Since $0 < \varepsilon < 1$, 
there exists a constant $c(t)$ independent of $x \in \mathbb R^N$ such that 
$u(t) = v(t) + c(t)$ in $\mathcal S'(\mathbb R^N)$. 
It follows from $\nabla u(t) = \nabla v(t)$ in $\mathcal S'(\mathbb R^N)$ that 
\begin{equation}\notag 
u(t) 
 = e^{t\mathcal L } u_0 
  + \int_0^t e^{(t-\tau) \mathcal L} |\nabla u(\tau)|^p\,d\tau 
 = e^{t\mathcal L } u_0 
  + \int_0^t e^{(t-\tau) \mathcal L} |\nabla v(\tau)|^p\,d\tau   
= v(t)
\quad \text{in }  \mathcal S'(\mathbb R^N). 
\end{equation}
Therefore, 
$c(t)\equiv0$ in $ \mathcal S'(\mathbb R^N)$,
and
the uniqueness follows. 
$\Box$\vspace{7pt}

\noindent 
{\bf Proof of the decay estimates \eqref{eq:1.9} and \eqref{eq:1.10}. }
According to the above proof of global existence, let 
\begin{equation}\
\label{eq:3.17}
\lambda = (p-1)/p,
\end{equation}
and it is sufficient to show only \eqref{eq:1.9} for the solution $u$ satisfying 
\begin{equation}\label{eq:3.18}
\| u \|_{\dot X^1_r \cap \dot X^1_\infty \cap \dot Y^\lambda_r \cap \dot Y^\lambda_\infty}  < \infty , 
\end{equation}
since \eqref{eq:1.10} is obtained by $\| u \|_{\dot Y^\lambda_r \cap \dot Y^\lambda _\infty} < \infty$. 
In the case $0 \leq t \leq 1$, the boundedness in time on 
the norms $\| \nabla^j u(t) \|_{L^q}$ ($j = 0,1$)
is obtained by the H\"older inequality, the inequalities  
$\| u(t) \|_{L^q} \preceq \| u(t) \|_{\dot B^0_{q,1}}$,
$\| \nabla u(t) \|_{L^q} \preceq \| u(t) \|_{\dot B^1_{q,1}}$, 
and \eqref{eq:3.18}, 
so that it suffices to consider the case $t > 1$. 
We show the estimate with derivative: 
\begin{equation}\label{eq:3.19}
\| \nabla u(t) \|_{L^q} 
\preceq (1+t) ^{- N(\frac{1}{r}-\frac{1}{q}) -1},\qquad t>0 . 
\end{equation}
Once \eqref{eq:3.19} is proved, 
it is possible to show the decay estimate of $\| u(t) \|_{L^q}$.
In fact,
by \eqref{eq:2.18} and \eqref{eq:3.19} we see that
\begin{equation}\notag 
\begin{split}
&
\| u(t) \|_{L^q} 
\\
& 
\preceq t^{-N(\frac{1}{r}-\frac{1}{q})}  \| u_0 \|_{L^r} 
   + \int_0^{t/2} (t-\tau)^{-N(\frac{1}{r}-\frac{1}{q})} \| |\nabla u (\tau)|^p \|_{L^r} \, d\tau
   + \int_{t/2}^t \| |\nabla u (\tau)|^p\|_{L^q} \, d\tau
\\
& 
\preceq t^{-N(\frac{1}{r}-\frac{1}{q})} 
   +  t^{-N(\frac{1}{r}-\frac{1}{q})}  
     \int_0^{t/2} \| \nabla u (\tau) \|_{L^{pr}}^p \, d\tau
   +  \int_{t/2}^t 
      \| \nabla u(\tau) \|_{L^{pq}} ^p \, d\tau 
\\
& 
\preceq t^{-N(\frac{1}{r}-\frac{1}{q})} 
   + t^{-N(\frac{1}{r}-\frac{1}{q})}  
     \int_0^{t/2} 
       \big\{ (1+\tau) ^{-N(\frac{1}{r}-\frac{1}{pr}) -1 } \big\} ^p 
     \,d\tau
   + \int_{t/2}^t
       \big\{ \tau ^{-N(\frac{1}{r}-\frac{1}{pq}) -1 } \big\} ^p
     \, d\tau 
\\
& 
\preceq t^{-N(\frac{1}{r}-\frac{1}{q})} ,\qquad t\ge1,
\end{split}
\end{equation}
so that, the decay estimate in $L^q (\mathbb R^n)$ is obtained.

We show \eqref{eq:3.19}. 
It follows from \eqref{eq:2.18} that 
\begin{equation}\label{eq:3.20}
\begin{split} 
\| \nabla u(t) \|_{L^ q} 
& 
\preceq t^{-N(\frac{1}{r}-\frac{1}{q})-1} \| u_0 \|_{L^r} 
    + \Big( \int_0^{t/2} + \int_{t/2}^t \Big) 
      \| \nabla e^{(t-\tau) \mathcal L} |\nabla u(\tau)|^p \|_{L^q}  \, d\tau.
\end{split}
\end{equation}
We first consider the case $r < \infty$. 
By \eqref{eq:2.18}, \eqref{eq:3.17} and \eqref{eq:3.18} we have 
\begin{equation}\label{eq:3.21}
\begin{split} 
  \int_0^{t/2} 
      \| \nabla e^{(t-\tau)\mathcal L} |\nabla u (\tau)|^p \|_{L^q} \, d\tau
& 
\preceq 
  \int_0^{t/2} 
     (t-\tau) ^{-N(\frac{1}{r} - \frac{1}{q}) -1}
      \| |\nabla u (\tau)|^p \|_{L^r} \, d\tau 
\\
& 
\preceq 
     t ^{-N(\frac{1}{r} - \frac{1}{q}) -1}
  \int_0^{t/2} 
      \| \nabla u (\tau) \|_{L^{pr} } ^p \, d\tau 
\\
& 
\preceq 
     t ^{-N(\frac{1}{r} - \frac{1}{q}) -1}
  \int_0^{t/2} 
    \big\{ (1+\tau)^{-N(\frac{1}{r} - \frac{1}{pr}) -1+\lambda} \big\}^p \, d\tau 
\\
& 
\preceq 
     t ^{-N(\frac{1}{r} - \frac{1}{q}) -1} . 
\end{split}
\end{equation}
On the other hand,
by \eqref{eq:3.17}, \eqref{eq:3.18},
the boundedness of $e^{(t-s)\mathcal L}$ in $L^q $
and the H\"older inequality 
we obtain 
\begin{equation}
\label{eq:3.22}
\begin{split} 
& 
  \int_{t/2} ^t
      \| \nabla e^{(t-\tau)\mathcal L} |\nabla u(\tau)|^p \|_{L^q} \,d\tau 
\\
& 
\preceq 
  \int_{t/2} ^t
      \| \nabla |\nabla u(\tau)|^p \|_{L^q} \,d\tau 
\\
& 
\preceq 
  \int_{t/2} ^t
      \| u(\tau) \|_{\dot B^1_{\infty,1}}^{p-1} \| u(\tau) \|_{\dot B^2_{q,1}} \,d\tau
\\
& 
\preceq 
  \| u \|_{\dot Y^\lambda_{\infty}} ^{p-1}
  \int_{t/2} ^t
     (\tau^{-\frac{N}{r} -1+\lambda})^{p-1} 
     \tau^{-N(\frac{1}{r}-\frac{1}{q}) -1+\lambda} 
     \tau^{N(\frac{1}{r}-\frac{1}{q}) +1-\lambda} 
     \| u(\tau) \|_{\dot B^2_{q,1}} \,d\tau 
\\
& 
\preceq 
     ( t^{-\frac{N}{r} -1+\lambda})^{p-1} 
     t^{-N(\frac{1}{r}-\frac{1}{q}) -1+\lambda} 
  \| u \|_{\dot Y^\lambda_{\infty}} ^{p-1}
  \| u \|_{\dot Y^\lambda_{q}}  
\\
& 
\preceq t^{-N(\frac{1}{r}-\frac{1}{q}) -1} . 
\end{split}
\end{equation}
This together with \eqref{eq:3.20} and \eqref{eq:3.21} yields \eqref{eq:3.19} for the case $r < \infty$.
 
Next we consider the case $r = \infty$.
In this case,
the problem is that the integral in the third line of \eqref{eq:3.21} diverges as $t \to \infty$. 
Then corresponding estimate to \eqref{eq:3.21} is the following 
with taking $r = q = \infty$
\begin{equation}\notag 
  \int_0^{t/2} 
      \| \nabla e^{(t-\tau)\mathcal L} |\nabla u (\tau)|^p \|_{L^\infty} \, d\tau 
\preceq 
     t ^{-1}
  \int_0^{t/2} 
    \big\{ (1+\tau)^{-1+\lambda} \big\}^p \, d\tau 
\preceq 
     t ^{-1} \log (1+t) ,
\end{equation}
and the same estimate as \eqref{eq:3.22} holds.
This implies that
$$
\| \nabla u(t) \|_{L^\infty} \preceq t^{-1} \log (1+t), \qquad t>1. 
$$
By this decay estimate, we can improve the corresponding one to \eqref{eq:3.21} as
\begin{equation}\notag 
  \int_0^{t/2} 
      \| \nabla e^{(t-\tau)\mathcal L} |\nabla u (\tau)|^p \|_{L^\infty} \, d\tau 
\preceq 
     t ^{-1}
  \int_0^{t/2} 
    \big\{ (1+\tau)^{-1} \log (2+\tau) \big\}^p \, d\tau 
\preceq 
     t ^{-1} . 
\end{equation}
Therefore, we also have the estimate \eqref{eq:3.19} for the case $r=\infty$,
and the proof of \eqref{eq:1.9} is completed. 
$\Box$\vspace{7pt}

\section{Asymptotic behavior}
In this section
we prove the assertion~(ii) of Theorem~\ref{Theorem:1.1}.
The proof is based on the arguments in \cite[Theorem~1.2]{IKK} and \cite[Theorem~1.1]{IKo}
(see, also, \cite{IK}).
Throughout this section we assume that
$u$ is a global-in-time solution of \eqref{eq:1.7} satisfying \eqref{eq:1.9} and \eqref{eq:1.10}.
\vspace{5pt}

\noindent
{\bf Proof of Theorem~\ref{Theorem:1.1}~(ii)-(a).}
Let $r\in(1,\infty)$ and $q\in[r,\infty]$.
By \eqref{eq:1.7}, for any $j \in \{0,1\}$, we have
\begin{equation}
\label{eq:4.1}
\begin{split}
&
\|\nabla^j[v(t)-e^{t\mathcal{L}}u_0]\|_{L^q}
\\
&
\le
\left\|\nabla^j\int_{t/2}^te^{(t-\tau){\mathcal L}}|\nabla v(\tau)|^p\,d\tau\right\|_{L^q}+\left\|\nabla^j\int_0^{t/2}e^{(t-\tau){\mathcal L}}|\nabla v(\tau)|^p\,d\tau\right\|_{L^q}
\\
&
=:I_{1,j}(t)+I_{2,j}(t)
\end{split}
\end{equation}
for all $t>0$.
We first estimate $I_{1,j}$.
By \eqref{eq:1.9} and \eqref{eq:2.18} we obtain
\begin{equation}
\label{eq:4.2}
\begin{split}
I_{1,0}(t)
&
\le\int_{t/2}^t\|e^{(t-\tau){\mathcal L}}|\nabla v(\tau)|^p\|_{L^q}\,d\tau
\\
&
\le\int_{t/2}^t\|\nabla v(\tau)\|_{L^\infty}^{p-1}\|\nabla v(\tau)\|_{L^q}\,d\tau
\\
&
\preceq\int_{t/2}^t\tau^{-(\frac{N}{r}+1)(p-1)}\tau^{-N(\frac{1}{r}-\frac{1}{q})-1}\,d\tau
\\
&
\preceq t^{-N(\frac{1}{r}-\frac{1}{q})-(\frac{N}{r}+1)(p-1)},\qquad t\ge1.
\end{split}
\end{equation}
Furthermore, applying the argument similar to \eqref{eq:3.22} with \eqref{eq:1.9} and \eqref{eq:1.10},
we see that
\begin{equation}
\label{eq:4.3}
\begin{split}
I_{1,1}(t)
&
\le\int_{t/2}^t\|\nabla e^{(t-\tau){\mathcal L}}|\nabla v(\tau)|^p\|_{L^q}\,d\tau
\\
&
\le\int_{t/2}^t\|\nabla v(\tau)\|_{L^\infty}^{p-1}\|\nabla^2 v(\tau)\|_{L^q}\,d\tau
\\
&
\preceq\int_{t/2}^t\tau^{-(\frac{N}{r}+1)(p-1)}\tau^{-N(\frac{1}{r}-\frac{1}{q})-\frac{1}{p}}\tau^{N(\frac{1}{r}-\frac{1}{q})+\frac{1}{p}}\|\nabla^2 v(\tau)\|_{L^q}\,d\tau
\\
&
\preceq t^{-N(\frac{1}{r}-\frac{1}{q})-\frac{1}{p}-(\frac{N}{r}+1)(p-1)}
\int_0^\infty \tau^{N(\frac{1}{r}-\frac{1}{q})+\frac{1}{p}}\|v(\tau)\|_{\dot B_{q,1}^2}\,d\tau
\\
&
\preceq t^{-N(\frac{1}{r}-\frac{1}{q})-\frac{1}{p}-(\frac{N}{r}+1)(p-1)},\qquad t\ge1.
\end{split}
\end{equation}
Since $(p-1)p+1>p$ for all $p>1$,
it follows from \eqref{eq:4.2} and \eqref{eq:4.3} that
\begin{equation}
\label{eq:4.4}
t^{N(\frac{1}{r}-\frac{1}{q})+j}I_{1,j}(t)=O(t^{-\frac{N}{r}(p-1)})
\end{equation}
as $t\to\infty$, for any $j\in \{0,1\}$.

Next we estimate $I_{2,j}$.
For the case $1<r\le p$, by \eqref{eq:1.9} and \eqref{eq:2.18} we obtain
\begin{equation}
\label{eq:4.5}
\begin{split}
I_{2,j}(t)
&
\le\int_0^{t/2}\|\nabla^je^{(t-\tau){\mathcal L}}|\nabla v(\tau)|^p\|_{L^q}\,d\tau
\\
&
\preceq\int_0^{t/2}(t-\tau)^{-N(1-\frac{1}{q})-j}\||\nabla v(\tau)|^p\|_{L^1}\,d\tau
\\
&
\preceq t^{-N(1-\frac{1}{q})-j}\int_0^{t/2}\|\nabla v(\tau)\|_{L^p}^p\,d\tau
\\
&
\preceq t^{-N(\frac{1}{r}-\frac{1}{q})-j-\frac{N}{r}(r-1)}\int_0^\infty(1+\tau)^{-N(\frac{p}{r}-1)-p}\,d\tau
\\
&
\preceq t^{-N(\frac{1}{r}-\frac{1}{q})-j-\frac{N}{r}(r-1)},\qquad t\ge1.
\end{split}
\end{equation}
For the case $r>p$, by \eqref{eq:1.9} and \eqref{eq:2.18} again we see that
\begin{equation*}
\begin{split}
I_{2,j}(t)
&
\le\int_0^{t/2}\|\nabla^je^{(t-\tau){\mathcal L}}|\nabla v(\tau)|^p\|_{L^q}\,d\tau
\\
&
\preceq\int_0^{t/2}(t-\tau)^{-N(\frac{p}{r}-\frac{1}{q})-j}\||\nabla v(\tau)|^p\|_{L^\frac{r}{p}}\,d\tau
\\
&
\preceq t^{-N(\frac{p}{r}-\frac{1}{q})-j}\int_0^{t/2}\|\nabla v(\tau)\|_{L^r}^p\,d\tau
\\
&
\preceq t^{-N(\frac{1}{r}-\frac{1}{q})-j-\frac{N}{r}(p-1)}\int_0^\infty(1+\tau)^{-p}\,d\tau
\preceq t^{-N(\frac{1}{r}-\frac{1}{q})-j-\frac{N}{r}(p-1)},\qquad t\ge1.
\end{split}
\end{equation*}
This together with \eqref{eq:4.5} yields
\begin{equation}
\label{eq:4.6}
t^{N(\frac{1}{r}-\frac{1}{q})-j}I_{2,j}(t)
=
\left\{
  \begin{array}{ll}
  O(t^{-\frac{N}{r}(r-1)})
  &
  \mbox{if}\quad p\ge r,
  \vspace{5pt}\\
  O(t^{-\frac{N}{r}(p-1)})
  &
  \mbox{if}\quad p<r,
  \end{array}
  \right.
\end{equation}
as $t\to\infty$, for any $j \in \{0,1\}$.
Therefore,
substituting \eqref{eq:4.4} and \eqref{eq:4.6} into \eqref{eq:4.1}, we have \eqref{eq:1.11},
and the assertion (ii)-(a) of Theorem~\ref{Theorem:1.1} follows.
$\Box$\vspace{7pt}
\noindent
{\bf Proof of Theorem~\ref{Theorem:1.1}~(ii)-(b).}
Let $r=1$.
Put
\begin{equation}
\label{eq:4.7}
c(t):=M(u_0) +\int_0^tM(|\nabla v(\tau)|^p)\,d\tau,
\end{equation}
where 
\begin{equation}
\label{eq:4.8}
M(f)=\int_{{\mathbb R}^N}f(x)\,dx.
\end{equation}
Then, by \eqref{eq:1.9} we have
$$
|c(t_2)-c(t_1)|= \int_{t_1}^{t_2}M(|\nabla v(\tau)|^p)\,d\tau
=\int_{t_1}^{t_2}\|\nabla v(\tau)\|_{L^p}^p\,d\tau
\preceq \int_{t_1}^{t_2}(1+\tau)^{-N(p-1)-p}\,d\tau
$$
for all $t_2\ge t_1\ge0$.
This implies that there exists the limit $C_*$ given by \eqref{eq:1.2} such that
\begin{equation}
\label{eq:4.9}
\left|\int_{{\mathbb R}^N}v(x,t)\,dx-C_*\right|=|c(t)-C_*|=O(t^{-(N+1)(p-1)})
\end{equation}
as $t\to\infty$.
Furthermore, \eqref{eq:2.17} and \eqref{eq:4.9} yield
\begin{equation}
\label{eq:4.10}
\lim_{t\to\infty}t^{N(1-\frac{1}{q})+j}\|\nabla^j[c(t)P_{t+1}-C_*P_{t+1}]\|_{L^q}
=\lim_{t\to\infty}|c(t)-C_*|=0
\end{equation}
for all $q\in[1,\infty]$ and $j=0,1$.
Put
\begin{equation}
\label{eq:4.11}
w(x,t):=[e^{t{\mathcal L}}u_0](x)-M(u_0)\,P_{t+1}(x)
\end{equation}
for all $(x,t)\in{\mathbb R}^N\times(0,\infty)$.
Since it follows from the semigroup property of $P_t$ that 
\begin{equation}
\label{eq:4.12}
[e^{t{\mathcal L}}P_1](x)=P_{t+1}(x),
\end{equation}
we have
$$
w(x,t)=[e^{t{\mathcal L}}w(0)](x).
$$
On the other hand, by \eqref{eq:1.5}, \eqref{eq:4.8} and \eqref{eq:4.11} we obtain
$$
\int_{{\mathbb R}^N}w(x,0)\,dx=0.
$$
Therefore, applying Lemma~\ref{Lemma:2.5} with the aid of \eqref{eq:2.18}, we see that
\begin{equation}
\label{eq:4.13}
\begin{split}
\lim_{t\to\infty}t^{N(1-\frac{1}{q})+j}\|\nabla^jw(t)\|_{L^q}
&
=\lim_{t\to\infty}t^{N(1-\frac{1}{q})+j}\|\nabla^je^{t{\mathcal L}}w(0)\|_{L^q}
\\
&
\preceq \lim_{t\to\infty}\|e^{\frac{t}{2}{\mathcal L}}w(0)\|_{L^1}=0
\end{split}
\end{equation}
for all $q\in[1,\infty]$ and $j=0,1$.

Let
\begin{equation}
\label{eq:4.14}
F(x,t):=|\nabla v(x,t)|^p-M(|\nabla v(t)|^p)\,P_{t+1}(x).
\end{equation}
Then, by \eqref{eq:1.5} and \eqref{eq:4.8} we have
\begin{equation}
\label{eq:4.15}
\int_{{\mathbb R^N}}F(x,t)\,dx=0,\qquad t\ge0.
\end{equation}
Since it follows from \eqref{eq:4.12} and \eqref{eq:4.14} that
\begin{equation*}
\begin{split}
&
\int_0^t e^{(t-\tau)\mathcal L}|\nabla v(\tau)|^p\,d\tau-\int_0^tM(|\nabla v(\tau)|^p)\,d\tau\,P_{t+1}(x)
\\
&
=\int_0^t e^{(t-\tau)\mathcal L}\left\{|\nabla v(\tau)|^p-M(|\nabla v(\tau)|^p)\,P_{\tau+1}\right\}\,d\tau
=\int_0^t e^{(t-\tau)\mathcal L}F(\tau)\,d\tau,
\end{split}
\end{equation*}
by \eqref{eq:1.7}, \eqref{eq:4.7} and \eqref{eq:4.11} we see that
\begin{equation*}
\begin{split}
&
v(x,t)-c(t)P_{t+1}(x)
\\
&
= e^{t{\mathcal L}}u_0+\int_0^t e^{(t-\tau)\mathcal L}|\nabla v(\tau)|^p\,d\tau-\left[M(u_0) +\int_0^tM(|\nabla v(\tau)|^p)\,d\tau\right]P_{t+1}(x)
\\
&
=w(x,t)+\int_0^t e^{(t-\tau)\mathcal L}F(\tau)\,d\tau.
\end{split}
\end{equation*}
This together with \eqref{eq:4.10} and \eqref{eq:4.13} implies that
\begin{equation*}
\begin{split}
&
\lim_{t\to\infty}t^{N(1-\frac{1}{q})+j}\|\nabla^j[v(t)-C_*P_{t+1}]\|_{L^q}
\\
&
=\lim_{t\to\infty}t^{N(1-\frac{1}{q})+j}\|\nabla^j[v(t)-c(t)P_{t+1}]\|_{L^q}
+\lim_{t\to\infty}t^{N(1-\frac{1}{q})+j}\|\nabla^j[c(t)P_{t+1}-C_*P_{t+1}]\|_{L^q}
\\
&
=\lim_{t\to\infty}t^{N(1-\frac{1}{q})+j}\left\|\nabla^j\int_0^t e^{(t-\tau)\mathcal L}F(\tau)\,d\tau\right\|_{L^q}.
\end{split}
\end{equation*}
Therefore, in order to obtain \eqref{eq:1.12}, it suffices to prove 
\begin{equation}
\label{eq:4.16}
\lim_{t\to\infty}t^{N(1-\frac{1}{q})+j}\left\|\nabla^j\int_0^t e^{(t-\tau)\mathcal L}F(\tau)\,d\tau\right\|_{L^q}=0. 
\end{equation}
For any $j\in \{0,1\}$, put
\begin{equation*}
\begin{split}
&
J_{1,j}(t):=\int_{t/2}^t\nabla^je^{(t-\tau){\mathcal L}}F(\tau)\,d\tau,
\\
&
J_{2,j}(t):=\int_L^{t/2}\nabla^je^{(t-\tau){\mathcal L}}F(\tau)\,d\tau,
\\
&
J_{3,j}(t):=\int_0^L\nabla^je^{(t-\tau){\mathcal L}}F(\tau)\,d\tau,
\end{split}
\end{equation*}
for $t\ge2L$, where $L\ge1$.
Since it follows from 
\eqref{eq:1.10}, \eqref{eq:2.17} and \eqref{eq:4.14} that
\begin{equation}
\label{eq:4.17}
\sup_{t>0}\,(1+t)^{N(1-\frac{1}{q})+N(p-1)+p}\|F(t)\|_{L^q}<\infty, 
\end{equation}
by \eqref{eq:2.18} we have
\begin{equation}
\label{eq:4.18}
t^{N(1-\frac{1}{q})}\|J_{1,0}(t)\|_q
\preceq t^{N(1-\frac{1}{q})}\int_{t/2}^t\|F(\tau)\|_{L^q}\,d\tau
\preceq 
t^{-(N+1)(p-1)}=o(1)
\end{equation}
as $t\to\infty$.
Furthermore, by \eqref{eq:1.9}, \eqref{eq:2.17} and \eqref{eq:4.4}
we obtain
\begin{equation}
\label{eq:4.19}
\begin{split}
t^{N(1-\frac{1}{q})+1}\|J_{1,1}(t)\|_{L^q}
&
\le
t^{N(1-\frac{1}{q})+1}\left[I_{1,1}(t)+\int_{t/2}^tM(|\nabla v(\tau)|^p)\|\nabla P_{\tau+1}\|_{L^q}\,d\tau\right]
\\
&
\preceq t^{N(1-\frac{1}{q})+1}I_{1,1}(t)+\int_{t/2}^t\tau^{-N(p-1)-p}\,d\tau
\\
&
\preceq t^{N(1-\frac{1}{q})+1}I_{1,1}(t)+t^{-(N+1)(p-1)}=o(1)
\end{split}
\end{equation}
as $t\to\infty$.
Moreover, by \eqref{eq:2.18} and \eqref{eq:4.17}
we have 
\begin{equation}
\label{eq:4.20}
\begin{split}
t^{N(1-\frac{1}{q})+j}\|J_{2,j}(t)\|_{L^q}
&
\le t^{N(1-\frac{1}{q})+j}\int_L^{t/2}(t-\tau)^{-N(1-\frac{1}{q})-j}\|F(\tau)\|_{L^1}\,d\tau\\
&
\preceq\int_L^{t/2}\|F(\tau)\|_{L^1}\,d\tau 
\preceq \int_L^{t/2}\tau^{-N(p-1)-p}\,d\tau
\preceq L^{-(N+1)(p-1)}
\end{split}
\end{equation}
for all sufficiently large $t$. 
Similarly, we see that 
\begin{equation}
\label{eq:4.21}
\begin{split}
 t^{N(1-\frac{1}{q})+j}\|J_{3,j}(t)\|_{L^q}
 &
 \le t^{N(1-\frac{1}{q})+j}\int_0^L\left\|\nabla^je^{\frac{(t-\tau)}{2}{\mathcal L}}e^{\frac{(t-\tau)}{2}{\mathcal L}}F(\tau)\right\|_{L^q}\,d\tau\\
 &
 \preceq\int_0^L\left\|e^{\frac{(t-\tau)}{2}{\mathcal L}}F(\tau)\right\|_{L^1}\,d\tau
\end{split}
\end{equation}
for all $t\ge 2L$. 
On the other hand, for any $L>0$, it follows from Lemma~\ref{Lemma:2.5} with \eqref{eq:4.15} that
\begin{equation}
\label{eq:4.22}
\lim_{t\to\infty}
\left\|e^{\frac{(t-\tau)}{2}{\mathcal L}}F(\tau)\right\|_{L^1}=0
\end{equation}
for all $s\in (0,L)$.
Furthermore, by \eqref{eq:2.18} we have
\begin{equation}
\label{eq:4.23}
\sup_{t\ge 2L}\left\|e^{\frac{(t-\tau)}{2}{\mathcal L}}F(\tau)\right\|_{L^1}
\le\|F(\tau)\|_{L^1}.
\end{equation}
Then, applying the Lebesgue dominated convergence theorem with \eqref{eq:4.22} and \eqref{eq:4.23} to \eqref{eq:4.21},
we obtain
\begin{equation}
\label{eq:4.24}
\lim_{t\to\infty}t^{N(1-\frac{1}{q})+j}\|J_{3,j}(t)\|_{L^q}=0.
\end{equation}
Therefore, by \eqref{eq:4.18}, \eqref{eq:4.19}, \eqref{eq:4.20} and \eqref{eq:4.24} 
we see that
$$
\limsup_{t\to\infty}t^{N(1-\frac{1}{q})+j}\left\|\nabla^j\int_0^t e^{(t-\tau)\mathcal L}F(\tau)\,d\tau\right\|_{L^q}
\le C_1L^{-(N+1)(p-1)}
$$
for some constant $C_1$ independent of $L$. 
Therefore, since $L$ is arbitrary, 
by $p>1$ we have \eqref{eq:4.16} ,
and the proof of the assertion (ii)-(b) of Theorem~\ref{Theorem:1.1} is complete.
$\Box$\vspace{7pt}

\noindent
{\bf Acknowledgements.}
The first author was supported by the Grant-in-Aid for Young Scientists (B) (No.~25800069)
from Japan Society for the Promotion of Science.
The second author was supported by the Grant-in-Aid for Young Scientists (B) (No.~24740107)
from Japan Society for the Promotion of Science.


\end{document}